\documentclass[a4paper,10pt]{omarams-l}

\copyrightinfo{2005}
  {M. Katsoulakis, G. Kossioris \& O. Lakkis}
\copyrightinfo{}{}

\PII{}
\usepackage{ifthen}
\provideboolean{isthesis}
\setboolean{isthesis}{false}
\provideboolean{issiamltex}
\setboolean{issiamltex}{false}
\provideboolean{isamsltex}
\setboolean{isamsltex}{true}
\provideboolean{pdfformat}%
\setboolean{pdfformat}{false}

\usepackage{amsmath,amssymb}    
\usepackage{stmaryrd} 
\usepackage{mathrsfs}
\newboolean{usemathrsfs}
\setboolean{usemathrsfs}{true}

\usepackage{enumerate}
\usepackage{xspace}
\usepackage{boxedminipage}
\usepackage[latin1]{inputenc}
\usepackage{verbatim}

\newcommand{\Ignore}[1]{}
\newcommand{\freeze}[1]{}
\provideboolean{shownotes}
\setboolean{shownotes}{true}
\newcommand{\margnote}[1]{
\ifthenelse{\boolean{shownotes}}{
\texttt{[*]}
\marginpar{\raggedright\tiny\texttt{[*] #1}}}{}}
\provideboolean{showchanges}
\setboolean{showchanges}{false}
\usepackage[dvips]{color}
\newcommand{\chcolor}{\color{blue}}
\newcommand{\changes}[1]{
  \ifthenelse{\boolean{showchanges}}
	     {{\chcolor #1}}
	     {#1}
}

\newcommand{\iid}{\ensuremath{\text{i.i.d.}}\xspace}

\usepackage{mathrsfs}
\provideboolean{usemathrsfs}
\setboolean{usemathrsfs}{true}

\ifthenelse{\boolean{usemathrsfs}}%
{\newcommand{\mathscript}{\mathscr}}
{\newcommand{\mathscript}{\mathcal}}
 
 \newcommand{\cB}{\ensuremath{\mathscript B}\xspace}
 
 \newcommand{\cD}{\ensuremath{\mathscript D}\xspace}
 
 \newcommand{\cF}{\ensuremath{\mathscript F}\xspace}

 \newcommand{\cI}{\ensuremath{\mathscript I}\xspace}

 \newcommand{\rE}{\ensuremath{\mathbb E}\xspace}
 \newcommand{\rN}{\ensuremath{\mathbb N}\xspace}
 
 \newcommand{\rR}{\ensuremath{\mathbb R}\xspace}

 \newcommand{\rV}{\ensuremath{\mathbb V}\xspace}

 \newcommand{\NO}{\ensuremath{\rN_0}}
 
 \newcommand{\reals}{\rR}

 \newcommand{\one}{\ensuremath{\mathbf 1}\xspace}

 \newcommand{\W}{\ensuremath{\Omega}\xspace}
 
 \newcommand{\epsi}{\ensuremath{\epsilon}\xspace}
 \newcommand{\ep}{\ensuremath{\varepsilon}\xspace}

 \newcommand{\lam}{\ensuremath{\lambda}\xspace}
 \newcommand{\ro}{\ensuremath{\varrho}\xspace}
 
 \newcommand{\fie}{\ensuremath{\varphi}\xspace}
 \newcommand{\w}{\ensuremath{\omega}\xspace}

 \newcommand{\norm}[1]{\ensuremath{\left|#1\right|}}
 \newcommand{\Norm}[1]{\ensuremath{\left\|#1\right\|}}
 
 \newcommand{\ltwop}[2]{\ensuremath{\left\langle#1,#2\right\rangle}}

 \newcommand{\crochet}[2]{\ensuremath{\left\langle #1,#2\right\rangle}}
 
 \newcommand{\ensemble}[2]{\ensuremath{\left\{ #1:\;#2 \right\}}}
 \newcommand{\setof}[1]{\ensuremath{\left\{#1\right\}}}
 
 \newcommand{\fromto}[2]{\ensuremath{\left[#1:#2\right]}}

 \renewcommand{\d}{\ensuremath{\,\mathrm{d}}}

 \newcommand{\registered}%
	    {\ensuremath{{}^{\bigcirc\!\;\!\!\!\!\!\!\!\;\text{\sc r}}}}

\newcommand{\lebmeas}[1]{\left|#1\right|}             

\newcommand{\Id}{\operatorname{Id}}                   

\newcommand{\Oh} {\operatorname{O}}                   
\newcommand{\pospart}[1]{\left[#1\right]^+}           

\renewcommand{\vec}[1]{\ensuremath{\boldsymbol #1}}
\newcommand{\prob}{\ensuremath{P}}

\newcommand{\Cov}{\ensuremath{\operatorname{Cov}}}

\newcommand{\var}{\ensuremath{\operatorname{Var}}}
\newcommand{\Var}{\ensuremath{\operatorname{Var}}}
\newcommand{\EX}{\operatorname{\rE}}
\newcommand{\ndistr}{\operatorname{N}}
\newcommand{\gaussrv}[2]{\ensuremath{\ndistr\left({#1},{#2}\right)}}

 \newcommand{\CC}{\ensuremath{\operatorname C}\xspace}
 \newcommand{\HH}{\ensuremath{\operatorname H}\xspace}
 \newcommand{\LL}{\ensuremath{\operatorname L}\xspace}

 \newcommand{\cont}[1]{\ensuremath{\CC^{#1}}}
 \newcommand{\contc}[1]{\ensuremath{\cont{#1}_{\mathrm{c}}}}

 \newcommand{\leb}[1]{\ensuremath{\LL_{#1}}}

 \newcommand{\sobh}[1]{\ensuremath{\HH^{#1}}}
 \newcommand{\sobhz}[1]{\ensuremath{\HH^{#1}_{0}}}

 \newcommand{\fes}[1]{\ensuremath{\rV^{#1}}}

\usepackage{amscd}

\newcommand{\imbedded}{\ensuremath{\hookrightarrow}}
\newcommand{\embeddedin}{\imbedded}

\newcommand{\secsymbol}{\S}
\newcommand{\secref}[1]{\secsymbol\ref{#1}}
\newcounter{LetterListItem}
\renewcommand{\theLetterListItem}{(\alph{LetterListItem})}
{
	\begin{list}%
	{\theLetterListItem\ }%
	{\usecounter{LetterListItem}
	 \ListParameters
	}
}%
{\end{list}}
\newcounter{NumberListItem}
\renewcommand{\theNumberListItem}{\arabic{NumberListItem}}
{
	\begin{list}%
	{\theNumberListItem.\ }%
	{\usecounter{NumberListItem}%
	 \ListParameters
	}
}%
{\end{list}}
\newcounter{StepsItem}

\newenvironment{Steps}%
{
	\begin{list}%
	{Step \theStepsItem.\ }%
	{\usecounter{StepsItem}%
	 \ListParameters
	}
}%
{\end{list}}
\providecommand{\grad}{\nabla}
\renewcommand{\grad}{\nabla}

\ifthenelse{\boolean{pdfformat}}
{
  \newcommand{\graphext}{pdf}
  \usepackage[pdftex]{graphicx} 
}
{
  \newcommand{\graphext}{eps}
  \usepackage[dvips]{graphicx} 
}

\usepackage{rotating}
\usepackage{amsthm} 
\ifthenelse{\boolean{isthesis}}{%
  \setcounter{secnumdepth}{1}%
}{
  \setcounter{secnumdepth}{2} 
}
\newcommand{\ListParameters}
{
	 \setlength{\topsep}{0em}
	 \setlength{\leftmargin}{0em}
         \setlength{\itemsep}{0ex}
	 \setlength{\parsep}{.5ex}
	 \setlength{\itemindent}{\labelsep}
	 \addtolength{\itemindent}{\labelwidth}
}
{\end{boxedminipage}}

\newtheoremstyle{plain}
  {}
  {}
  {\mdseries\slshape}
  {}
  {\bfseries}
  {}
  {.5ex}
  {}

\newtheoremstyle{note}
  {}
  {}
  {}
  {}
  {\bfseries}
  {}
  {.5ex}
  {}

\newtheoremstyle{claim}
  {}
  {}
  {\mdseries\slshape}
  {\parindent}
  {\bfseries}
  {}
  {.5ex}
  {}

\newtheoremstyle{exercise}
  {}
  {}
  {}
  {}
  {\bfseries}
  {}
  {1em}
  {}

  \newcommand{\ObsName}{Remark}
  \newcommand{\DefName}{Definition}
  \newcommand{\ExaName}{Example}
  
  \newcommand{\TheName}{Theorem}
  \newcommand{\LemName}{Lemma}
  \newcommand{\ProName}{Proposition}
  \newcommand{\CorName}{Corollary}
  \newcommand{\PbmName}{Problem}
  \newcommand{\HypName}{Hypothesis}
  
  \newcommand{\ExeName}{Exercise}
  \newcommand{\SolName}{Solution of problem}%
  \newcommand{\ClaName}{Claim}
  
  \newcommand{\Proofname}{Proof}

\ifthenelse{\boolean{isthesis}}{%
  \newcommand{\Thecounter}{The}
}{%
  \newcommand{\Thecounter}{subsection}
}

{\swapnumbers
  {
    \theoremstyle{plain}
    \ifthenelse{\boolean{isthesis}}{%
      \newtheorem{The}{\TheName}[section]
    }{%
      \newtheorem{The}[subsection]{\TheName}
    }
    
    {
      \theoremstyle{plain}

      \newtheorem{Lem}[\Thecounter]{\LemName}

    }
    
    {
      \theoremstyle{note}
      \newtheorem{Obs}[\Thecounter]{\ObsName}
      \newtheorem{Def}[\Thecounter]{\DefName}

    }
  }
}

{\theoremstyle{exercise}
\newtheorem{Exe}{\ExeName}[section]

}

{\begin{Exe}\vspace{-24pt}\begin{enumerate}[(a)\ ]}%
{\end{enumerate}\end{Exe}}

\renewcommand{\qed}{\vrule height 5pt depth 0pt width 3pt}

\newcounter{passo}

\newenvironment{Proof}[1][{}]%
{\noindent{\bf \Proofname\ #1}\setcounter{passo}{0}}%
{{\raggedright{{ }\hfill\qed}}} 

\newcommand{\fep}{\ensuremath{f_\epsi}}
\newcommand{\sep}{\ensuremath{\epsi^\gamma}}
\newcommand{\Fep}{\ensuremath{F_\epsi}}

\newcommand{\baru}{\ensuremath{\bar u}\xspace}
\newcommand{\bare}{\ensuremath{\bar e}\xspace}
\newcommand{\barr}{\ensuremath{\bar r}\xspace}
\newcommand{\ut}{\ensuremath{\partial_t u\xspace}}

\newcommand{\baret}{\ensuremath{\partial_t\bare}\xspace}
\newcommand{\lap}{\ensuremath{\Delta}}
\providecommand{\bm}{\ensuremath{W}\xspace}
\renewcommand{\bm}{\ensuremath{W}\xspace}
\newcommand{\hatbm}{\ensuremath{\widehat{\bm}}\xspace}
\newcommand{\barbm}{\ensuremath{\bar{\bm}}\xspace}
\newcommand{\noise}{\ensuremath{\partial_{xt}\bm\xspace}}
\newcommand{\barnoise}{\ensuremath{\partial_{xt}\barbm}}
\newcommand{\hatnoise}{\ensuremath{\partial_{xt}\hatbm}}
\newcommand{\bareta}{\ensuremath{\bar\eta}}

\newcommand{\etamn}{\ensuremath{\eta_{m,n}}}
\newcommand{\baretamn}{\ensuremath{\bareta_{m,n}}}

\newcommand{\abil}[2]{\ensuremath{\crochet{\grad #1}{\grad #2}}}
\newcommand{\xm}{\ensuremath{x_m}}
\newcommand{\xmo}{\ensuremath{x_{m-1}}}
\newcommand{\tn}{\ensuremath{{t_n}}}

\newcommand{\tno}{\ensuremath{{t_{n-1}}}}

\newcommand{\tparti}{\ensuremath{\cI_\rho}}
\newcommand{\sparti}{\ensuremath{\cD_\sigma}}
\newcommand{\stparti}{\ensuremath{\sparti\times\tparti}}
\newcommand{\thparti}{\ensuremath{\cI_k}}
\newcommand{\shparti}{\ensuremath{\cD_h}}

\newcommand{\honezd}{\ensuremath{\sobh 1_0(D)}}
\newcommand{\awn}{AWN\xspace}

\newcommand{\ito}{It\^o\xspace}
\renewcommand{\grad}{{\partial_x}}
\renewcommand{\lap}{{\partial_{xx}}}

\usepackage{algorithmic}
\usepackage{subfigure}
\numberwithin{equation}{section}
\begin{document}
%
\newcommand{\mytitle}%
{NOISE REGULARIZATION AND COMPUTATIONS FOR THE $1$-DIMENSIONAL STOCHASTIC ALLEN--CAHN PROBLEM}
\ifthenelse{\boolean{isamsltex}}%
{ \title[Regularization and computations for the stochastic Allen--Cahn]%
{\mytitle}
\author{Markos A. Katsoulakis}
\address{ %
  Markos A. Katsoulakis\newline
  Department of Mathematics and Statistics\newline
  University of Massachusetts\newline
  Amherst, MA 01003, USA
}
\curraddr{}
\email{markos@math.umass.edu}
\thanks {M.A.K.'s research is partially supported by the
NSF grants DMS-ITR-0219211  and DMS-0413864}
\author{Georgios T. Kossioris}
\address{ %
  Georgios T. Kossioris\newline
  Department of Mathematics\newline %
  University of Crete\newline %
  Knossos Avenue\newline
  Heraklion GR-71409, Greece\newline
  and\newline
  Institute for Applied and Computational Mathematics\newline
  Foundation for Research and Technology-Hellas\newline
  Vasilika Vouton\newline
  Crete GR-71110, Greece
}
\curraddr{}
\email{kosioris@math.uoc.gr}
\thanks{G.T.K.'s research is partially supported by 
  the EU's \emph{Research Training Network 
  on Fronts and Singularities} HPRN-CT-2002-00274.}
\author{Omar Lakkis} 
\address{ Omar Lakkis\newline
  Department of Mathematics\newline
  University of Sussex\newline
  Brighton, UK-BN1 9RF, United Kingdom} 
\curraddr{}
\email{o.lakkis@sussex.ac.uk} 
\thanks{
  O.L. was supported by the EU's \emph{Research Training Network on
  Hyperbolic and Kinetic Equations (HYKE)} HPRN-CT-2002-00282 and the
  EU's \emph{MCWave Marie Curie Fellowship} HPMD-CT-2001-00121 during
  his stay in Heraklion, and by The Nuffield Foundation's young
  researcher grant in Brighton.
}
\subjclass[2000]{}
\date{August 2005}
\commby{}
\dedicatory{}
}
{
\title{\mytitle}
\author{M.~A.~Katsoulakis, G.~T.~Kossioris and O.~Lakkis}
\date{\today}
}
\maketitle
\begin{abstract}
We address the numerical discretization of the Allen-Cahn problem with
additive white noise in one-dimensional space.  Our main focus is to
understand the behavior of the discretized equation with respect to a
small ``interface thickness'' parameter and the noise intensity.  The
discretization is conducted in two stages: (1) regularize the white
noise and study the regularized problem, (2) approximate the
regularized problem.  We address (1) by introducing a piecewise
constant random approximation of the white noise with respect to a
space-time mesh.  We analyze the regularized problem and study its
relation to both the original problem and the deterministic Allen-Cahn
problem.  Step (2) is then performed leading to a practical
Monte-Carlo method combined with a Finite Element-Implicit Euler
scheme.  The resulting numerical scheme is tested against theoretical
benchmark results concerning the behavior of the solution as the
interface thickness goes to zero.
\end{abstract}
%
%
\section{Introduction}
Stochastic partial differential equation (SPDE) models arise in
numerous applications ranging from materials science, surface
processes and macromolecular dynamics \cite{cook,spohn}, to atmosphere
and ocean modeling \cite{lin-neelin:03} and epidemiology
\cite{durrett}. These models are typically derived from finer and more
detailed models where unresolved degrees of freedom are represented by
suitable stochastic forcing terms. There are also some notable
rigorous derivations from microscopic scales in special asymptotic
regimes \cite[e.g.]{presutti,tribe}.

An important class of models consists of the stochastic
Ginzburg--Landau models which are typically obtained from microscopic
lattice models for a suitable order parameter (e.g., spin), by
statistical mechanics renormalization arguments combined with detailed
balance laws. 

Numerical simulation of these nonlinear SPDE's constitutes an
important research issue.  On the practical side, one is interested in
having efficient, reliable and not too complex numerical codes which
can be used either in the context of Monte Carlo methods or for sample
paths simulations of the physical models of phase transition
\cite[e.g.]{warren-boettinger:95,kielhorn-muthukumar:99,katsoulakis-kho:01,shardlow:00}.
From a more theoretical view-point, understanding the issues arising
from the discretization of SPDE's, in a more general setting than
phase separation, both through finite difference or finite element
schemes, turns out to be a non-obvious departure from numerical
schemes for deterministic models
\cite{allen-novosel-zhang,gyongy:99,du-zhang:02,babuska-tempone-zouraris:04,schwab-todor:03}.

In this paper, our focus is on the numerical simulation of the
stochastic Allen--Cahn problem, which is one of the simplest models
exhibiting the phenomena of interface formation and nucleation.  The
stochastic Allen--Cahn problem is an ad-hoc white noise perturbation
of the deterministic Allen--Cahn, given by
\begin{equation}
  \label{eqn:noisy-allen-cahn}
  \ut(x,t)-\lap u(x,t)+\fep(u(x,t))=\sep\noise(x,t),
  \text{ for }x\in D,t\in[0,\infty)
\end{equation}
where $D=(-1,1)\subset\reals^1$, $\epsi>0$ and $\fep$ is an odd
nonlinearity scaled by $1/\epsi^2$ and $\noise$ is the space-time
white noise (see \secref{sse:setup.sac.problem} the details).  This is a stochastic version of the well-known deterministic
Allen--Cahn problem describing the evolution in time of a
polycrystalline material \cite{allen-cahn}.  We take boundary
conditions of Neumann type and the initial condition to be a
\emph{resolved profile}; we refer again to \secref{sec:setup} for the
details.  Note that this equation, with white noise, is tractable only
in $1$ spatial dimension, which is the case we will study.  In higher
space dimension, one has to consider noise which is colored in
space.

Equation \eqref{eqn:noisy-allen-cahn} is a type A model in Halperin's
classification \cite{hohenberg-halperin:77}.  It is non-conservative
in the order parameter $u$ and exhibits both nucleation and interface
formation, whilst retaining a relatively simple structure without
multiplicative or conservative noise terms encountered in type B models, such as the
Cahn--Hilliard--Cook equation \cite{kielhorn-muthukumar:99}.
 
While a thorough discussion of \eqref{eqn:noisy-allen-cahn} is given
in \secref{sec:setup}, it is worth mentioning here that this SPDE, with the
white noise term, is well-posed only in one space dimension.  Two
important pieces of work concerned with the analytic and probabilistic
aspects of \eqref{eqn:noisy-allen-cahn} are those of Funaki
\cite{funaki:95} and Brassesco, De Masi \& Presutti
\cite{brassesco-demasi-presutti:95}.  In both papers, the authors
study the asymptotic behavior of the solution processes as
$\epsilon\rightarrow0$.  In particular, it turns out that, under
suitable time-space rescaling, the solution with initial value taken
to be (roughly speaking) a step function, converges (in an appropriate
probabilistic sense) to the step function with its jump point
performing a Brownian motion. 

Though finite difference schemes have been used for simulations
\cite{kielhorn-muthukumar:99,katsoulakis-kho:01}, we follow here a
finite element approach.  The reason driving us to understand finite
element methods (FEM) for such equations is that FEM constitute a
quite flexible tool, especially for problems in higher dimensions
where one may have to deal complex geometries.  Also, finite elements
are naturally suited for adaptive schemes where fine scales may be
resolved only on small portions of the domain in order to obtain a
reasonable accuracy.  We believe that understanding the FEM in a
non-adaptive one-dimensional setting will pave the way to more
sophisticated studies.

Our strategy to formulate a finite element scheme for
\eqref{eqn:noisy-allen-cahn}, follows an idea introduced for linear
problems by Allen, Novosel \& Zhang \cite{allen-novosel-zhang}, and
consists in two steps:
\begin{enumerate}[1.\ ]
\item \emph{regularize the noise} term $\noise$, by replacing it with
  a somewhat smoother approximate white noise $\barnoise$;
\item \emph{discretize the regularized} problem. 
\end{enumerate}
This approach allows us to conduct a rigorous analysis of the
approximation.  It makes the subsequent finite element discretization
straightforward.  Note that a finite difference variant based on our
regularization is also possible.

Our first task, carried out in \secref{sec:approximation}, is to
construct a regularization, denoted $\barnoise(x,t)$, of $\noise(x,t)$
(appearing in \eqref{eqn:noisy-allen-cahn}) with respect to an
underlying uniform partition, \stparti, of the space-time domain
$D\times I$.  In the spirit of FEM, this regularization process
consists of a projection of the white noise onto an appropriate space
of piecewise constant space-time functions, which may be viewed as the
mixed derivatives of hat functions.  This idea, which has been
successfully used in the context of the linear heat equation
\cite{allen-novosel-zhang}, leads to the \emph{regularized problem}
\begin{equation}
  \label{eqn:noisy-allen-cahn-regularized}
  \ut(x,t)-\lap u(x,t)+\fep(u(x,t))=\sep\barnoise(x,t),
  \text{ for }x\in D,t\in[0,\infty).
\end{equation}
Notice that $\barnoise$ is still a stochastic process in space-time,
but it is much smoother than the white noise which allows equation
\eqref{eqn:noisy-allen-cahn-regularized} to be interpreted in the
usual PDE sense pathwise.  

In \secref{sec:convergence}, after recalling some basic properties of
problem \eqref{eqn:noisy-allen-cahn-regularized} and its solution, we
prove Theorem \ref{the:awn-exact-convergence}, which states that the
solution of the regularized problem converges---in an appropriate
sense---to the solution of the original SPDE
\eqref{eqn:noisy-allen-cahn} as the space-time partition becomes
infinitely fine.

Next, in \secref{sec:weak.noise}, we relate the solution of the
regularized problem to the deterministic solution of the Allen--Cahn
equation.  Our main result here, Theorem \ref{the:low.noise}, proved
for $\gamma>3$, indicates that the regularization parameters have to
be sufficiently small for the noise to be captured in the numerical
computations.  In fact, according to this Theorem the weaker the
noise, the finer one must take the space-time mesh, in order to see
the noise effects.  This is due to the fact that for a fixed
space-time mesh and $\epsi\to0$, the distance between the regularized
stochastic solution $\baru$ and the deterministic solution, $q$, is of
higher order in $\epsi$ than the distance between $\baru$ and $u$.  Our proof
makes use of the \emph{spectrum estimates} of the linearized elliptic
differential operator $-\lap+\fep'(q)$, derived independently by Xinfu
Chen \cite{chen.xinfu:94} and de Mottoni \& Schatzman
\cite{demottoni-schatzman:95}

We note that while numerical schemes for the stochastic Allen--Cahn
involving a spectral approach to white noise have been analyzed
\cite{liu.di:03} this is, up to our knowledge, a first analysis using
projection methods to regularize the white noise.

Step 2 of our strategy is accomplished in \secref{sec:fem}, where we
derive a simple finite element scheme for the regularized problem
\eqref{eqn:noisy-allen-cahn-regularized}.  This is a scheme which uses
piecewise polynomial finite elements to discretize the space variable
and an implicit (backward) Euler scheme to discretize the time
variable. Related numerical schemes have been thoroughly analyzed and
successfully applied in the context of the \emph{deterministic
Allen--Cahn problem}
\cite{feng-prohl:03,kessler-nochetto-schmidt:03,feng-wu:05} and for
the \emph{stochastic linear heat diffusion problem}
\cite{allen-novosel-zhang}.  It is for the first time, up to our
knowledge, that this scheme is employed in a \emph{stochastic and
nonlinear} setting.  The issues of regularity of the regularized
solution and the convergence of the FEM are objects of our current
research.

In \secref{sec:numerics}, we test our scheme in combination with a
Monte Carlo simulation.  The test consists in reconciling the
computational results with the theoretical results obtained by Funaki
\cite{funaki:95} and Brassesco, De Masi \& Presutti
\cite{brassesco-demasi-presutti:95} independently.  Our benchmarking
procedure consists in tracking the so-called center of a resolved
profile of the Allen--Cahn equation as time evolves, performing
statistics thereon and comparing them with the probabilistic results
coming from the theory.  The following conclusions are drawn: (1) The
robustness of the Monte Carlo method depends on the noise intensity,
the lower is the noise the higher is the observed robustness. (2) The
noise has to be resolved satisfactorily in order to see stochastic
effects. In contrast with the first conclusion, the lower the noise the
more has one to resolve the mesh in order to see the noise.  This is
in competition with the need to have a fine mesh in order to resolve
the transition layer, due to the structure of the solution of the
Allen--Cahn equation.  (3) The behavior captured by the numerics is
consistent with the theoretical results; in particular, the
Mueller--Funaki time scale $1+2\gamma$ (see
\ref{sse:low.noise.resolution} for the details) and the corresponding
Brownian motion diffusion coefficient are clearly exhibited by our
numerical results.  We close with some computations that capture the
drift of the interface, modeled by the Allen--Cahn equation.  This
drift, typical of the stochastic solution is quite fast with respect
to the deterministic case where the solutions are metastable states.

\section{Set up}
\label{sec:setup}
\subsection{Noisy Allen--Cahn problem}
\label{sse:setup.sac.problem}
We will study an initial-boundary value problem associated with the
semilinear parabolic partial differential equation with additive white
noise, known as the \emph{stochastic} (or \emph{noisy})
\emph{Allen--Cahn equation} given by \eqref{eqn:noisy-allen-cahn}.  The
nonlinearity $\fep$ is the derivative of an even coercive
function $\Fep$ with exactly two minimum points.  A function such as $\Fep$ is
known as a \emph{double-well potential} and, for sake of conciseness,
we focus on the model potential explicitly defined by
\begin{equation}
  \Fep(\xi)=\frac1{4\epsi^2}(\xi^2-1)^2,\text{ for }\xi\in\reals.
\end{equation}
Here $\epsi\in\reals^+$ is a scaling parameter.  The term $\noise$ is
the space-time {\em Gaussian white noise}, which can be defined as the
mixed distributional derivative of a {\em Brownian sheet} \bm
\cite{walsh:notes:84,kallianpur-xiong:notes:95}. The parameter
$\gamma\in\reals$ models the intensity of the white noise and plays a
delicate role in the analysis, as $\epsi\rightarrow0$.

The presence of the right-hand side makes \eqref{eqn:noisy-allen-cahn}
a randomly perturbed version of the Allen--Cahn equation which is a
{\em stochastic PDE} (SPDE).  A solution of such an equation has to be
interpreted in the stochastic sense.  That is, for each $t$, the
solution $u(\cdot,t)$ is understood as a random process on an
underlying probability measure space $(\W,\cF,\prob)$ with values in
a suitable function space defined on $D$.  Equation
\eqref{eqn:noisy-allen-cahn}, supplemented with the initial condition
\begin{equation}
  \label{eqn:initial-condition}
  u(x,0)=u_0(x),\:\forall x\in D,
\end{equation}
and with the Neumann boundary conditions
\begin{equation}
  \label{eqn:boundary-conditions}
  \partial_x u(-1,t)=\partial_x u(1,t)=0\quad\forall t\in \reals^+,
\end{equation}
defines the \emph{stochastic Allen--Cahn problem}.  For simplicity, we
assume that the initial condition $u_0$ is smooth enough and satisfies
the boundary conditions.  In \secref{sec:weak.noise} we shall focus on
a more particular class of initial conditions known as resolved profiles.
\subsection{Space-time stochastic integral}
One can give a mathematically rigorous definition of a solution of the
stochastic Allen--Cahn problem
\eqref{eqn:noisy-allen-cahn},\eqref{eqn:initial-condition}--\eqref{eqn:boundary-conditions}
as a distribution-valued process
\cite{walsh:notes:84,kallianpur-xiong:notes:95}.  However, we find it
more convenient, as in the case of the white noise generated from a
Brownian motion, to work with the \emph{stochastic integral} with
respect to the Brownian sheet $W$ denoted by ``$\int\cdot\d\bm$''
\cite[\S II]{walsh:notes:84} \cite[Ch. 3]{kallianpur-xiong:notes:95}.
In our doing so, we bear in mind the formal relationship
\begin{equation}
  \label{eqn:formal-integral}
  \int_0^\infty\int_D f(x,t)\noise(x,t)\d x\d t
  =
  \int_0^\infty\int_D f(x,t)\d\bm(x,t)
\end{equation}
that will inspire the weak formulation \eqref{eqn:weak.formulation}
and the definitions in \secref{sec:approximation}.  In the particular
case where $f$ is the characteristic function of a Borel-measurable
set $A \in \cB(\reals^+ \times D)$ of Lebesgue measure $\lebmeas
A<\infty$ the following basic property of the stochastic integral is
satisfied
\begin{equation}
  \label{eqn:white.noise.is.gaussian}
  \int_0^\infty\int_A \d\bm(x,t)=W(A)\in\gaussrv0{\lebmeas{A}},
\end{equation}
i.e., $W(A)$ is a Gaussian random variable with mean zero and variance
$\lebmeas{A}$.\footnote{For $\mu\in\reals$,
$\sigma\in\reals^+$ we denote by $\gaussrv\mu{\sigma^2}$ the class of
normally distributed (or Gaussian) random variables of mean $\mu$ and
variance $\sigma^2$ on the space $\W$.}

Since we are interested in numerical solutions, we consider the time
domain to be a bounded interval $I=[0,T]$, for some fixed
$T\in\reals^+$.  A fundamental property of the stochastic integral is
the following well-known \emph{$\leb2$-isometry}, which holds for the \ito
integral,
\begin{equation}
  \label{eqn:ito.isometry}
  \EX\left[\left(\int_I\int_D f(x,t)\d\bm(x,t)\right)^2\right]
  =\EX\left[\int_I\int_D f(x,t)^2\d x\d t\right],
\end{equation}
for any $\cF_t^W$-measurable $f\in\leb2(I\times D\times\W)$, where 
\begin{equation}
  \cF_t^W=\sigma\ensemble{W(A)}{A \in \cB(I \times D)},
\end{equation}
is the sigma-field (or sigma-algebra) generated by $\bm$ up to time
$t$, and $\EX$ denotes the \emph{expectation} with respect to
$(\W,\cF,\prob)$.\footnote{In compliance with the standard practice in
stochastic differential equations, we write explicitly the
probability variable $\w\in\W$ as an argument to random variables only
when necessary in order to avoid confusion.}

A useful consequence of \eqref{eqn:ito.isometry} is that
\begin{equation}
  \label{eqn:product.formula}
  \EX\left[\int_I\int_D f(x,t)\d\bm(x,t)\int_I\int_D
  g(y,s)\d\bm(y,s)\right]
  =\EX\left[\int_I\int_D f(x,t)g(x,t)\d x \d t\right],
\end{equation}
for any $\cF_t^W$-measurable $f,g\in\leb2(I\times D\times\W)$. In the
special case where $f$ and $g$ are, respectively, the characteristic
functions of two Borel sets $A$ and $B \in \cB(I \times D)$, with
$\lebmeas A,\lebmeas B<\infty$, \eqref{eqn:product.formula} implies
\begin{equation}
  \Cov(W(A),W(B))=\lebmeas{A\cap B}.
\end{equation}
\subsection{Integral solutions}

By multiplying \eqref{eqn:noisy-allen-cahn} with a test function
$\phi\in \contc2(D\times (0,\infty))$ and using the formal relation
\eqref{eqn:formal-integral}, one can write the problem in the usual
weak form\footnote{Whenever the meaning is clear from the context, for
sake of conciseness, we often drop the variables ``$x,t$'' and, in
non-stochastic integrals, also the corresponding elementary terms
``$\d$''.}
\begin{equation}
  \label{eqn:weak.formulation}
  \int_0^\infty
  \int_D \left(u \partial_t\phi
  - \grad u \grad \phi
  - \fep(u)\phi\right)
  +\sep \int_0^\infty\int_D \phi \d W=0.
\end{equation}
Despite the above formulation being quite useful, especially for
studying a numerical scheme, it is not very convenient to nail down
the concept of solution.  A rather more convenient way to give
rigorous meaning to \eqref{eqn:noisy-allen-cahn} is to look for an
integral solution of an \emph{equivalent integral equation}
\cite{daprato-zabczyk:book:92,doering:87,faris-jonalasinio:82,walsh:notes:84},
as we briefly illustrate next.

Introduce first the corresponding boundary value problem for the
stochastic linear heat equation
\cite{daprato-zabczyk:book:92,walsh:notes:84}
\begin{equation}
  \begin{gathered}
    \partial_t Z-\lap Z=\noise,\text{ in }D\times\reals^+_0\\
    Z(x,0)=0,\text{ on }D\\
    \grad Z(1,t)=\grad Z(-1,t)=0,\:\forall t\in[0,\infty).
  \end{gathered}
\end{equation}
The solution to this problem can be defined as the Gaussian process in
space-time produced by the stochastic integral
\begin{equation}
  Z_t(x)=Z(x,t):=\int_0^t\int_D G_{t-s}(x,y)\d\bm(y,s),
\end{equation}
where $G$ is the heat kernel for the corresponding homogeneous Neumann
problem.  In our one-dimensional particular case, $G$ can be
explicitly written as
\begin{equation}
  \label{eqn:heat.kernel}
  G_t(x,y)
  =4\sum_{k=0}^\infty (2-\delta_0^k)
  \cos\frac{\pi k(x+1)}2
  \cos\frac{\pi k(y+1)}2
  \exp\frac{-\pi^2 k^2 t}4,
\end{equation}
where $\delta_0^k$ is the Kronecker symbol.

 The \emph{integral solution} of \eqref{eqn:noisy-allen-cahn} can then
 be defined as a solution of the equivalent integral equation
\begin{equation}
  \label{eqn:integral.noisy.allen-cahn}
  \begin{split}
    u(x,t)=-\int_0^t\int_D G_{t-s}(x,y)\fep(u(y,s))\d y\d s
    +\int_D G_t(x,y)u_0(y)\d y
    +\sep Z_t(x).
  \end{split}
\end{equation}
It is known that such a solution exists uniquely as a
$\cont0(D)$-valued continuous process, $t\mapsto u(.,t)$, adapted to
$Z_t$, provided the initial condition $u_0$ satisfies the Neumann
boundary conditions
\cite{brassesco-demasi-presutti:95,walsh:notes:84,faris-jonalasinio:82}.
In this article we use this concept of solution which we refer to
simply as the \emph{solution of Problem
\eqref{eqn:noisy-allen-cahn},\eqref{eqn:initial-condition}--\eqref{eqn:boundary-conditions}}
and we will denote it by $u$.  Notice that $u$ is also referred to by
some authors as the \emph{Ginzburg-Landau process}
\cite{brassesco-demasi-presutti:95}.


For the aims set in this paper, namely, in order to study the error of
convergence of an approximation of the solution of
\eqref{eqn:noisy-allen-cahn}, we will need a uniform bound for $u$.
While in the deterministic case such a bound is direct consequence of
the maximum principle, in the stochastic case one cannot expect to
have a uniform bound in the whole probability space. However, a bound
on a set with large probability controlled by $\epsi$ will suffice for
our needs.  We present an extension of a previously known result
of Brassesco et al. \cite[Pro.5.2]{brassesco-demasi-presutti:95}.
\begin{Lem}[Probabilistic maximum principle]
  \label{lem:prob.max.princ}
  Let $\gamma>-1/2$.
  For each $T>0$ and $K_0>0$ there exist $c_1,\,c_2,\delta_0>0$ such that if 
  $\Norm{u_0}_{\leb\infty(D)}\leq 1+\delta_0$ then
  \begin{equation}
    \prob\setof{\sup_{t\in[0,T]}\Norm{u(t)}_{\leb\infty(D)}>1+K_0}
    \leq c_1\exp(-c_2/\epsi^{1+2\gamma}).
  \end{equation}
\end{Lem}
\begin{Proof}
  We reduce the proof to that of
  \cite[Pro.5.2]{brassesco-demasi-presutti:95} by introducing the
  time-space rescaling: $t\mapsto t/\epsi^2$ and $x\mapsto
  x/(\sqrt2\epsi)$ and extending the solution periodically to the
  whole space as to obtain the proper barrier function.  Since we are
  dealing with the more general case $\gamma>-1/2$, while they deal
  with the case $\gamma=0$ only, we retrace the salient points of
  their proof.  The barrier function $v$ satisfies the following
  equation---corresponding to
  \cite[(5.12)]{brassesco-demasi-presutti:95}:
  \begin{equation}
    \partial_t v -\frac12\lap v+2v
    =
    -3v^2-v^3+2^{-1/4}\epsi^{\gamma+1/2}\noise.
  \end{equation}
  Consider now the function 
  \begin{equation}
    V(x,t)=\int_0^t\exp(-2(t-s))H_{t-s}^{(\epsi\sqrt2)}(x,y)\d\bm(y,s),
  \end{equation}
  where $H_{t-s}^{(\epsi\sqrt2)}$ is  the Green operator defined by 
  \begin{equation}
    \exp(-2t)H_{t-s}^{(\epsi\sqrt2)}
    =\left(\partial_t -\frac12\lap +2\Id \right)^{-1},
  \end{equation}
  with homogeneous boundary conditions on
  $\left(-1/(\sqrt2\epsi),1/(\sqrt2\epsi)\right)$.  By using equation
  \cite[(5.2)]{brassesco-demasi-presutti:95} with $\lam
  =\exp(-(\gamma+1/2))$ and adapting properly the proof of \cite[Lemma
  2.1]{brassesco-demasi-presutti:95} we can easily conclude that for
  each $b>0$ there exist $c_1$ and $c_2>0$ such that
  \begin{equation}
    \prob^\epsi\setof{\sup_{t\leq T\epsi^{-2},x\in\reals}
    \norm{\epsi^{\gamma+1/2}V(x,t)}>b}
    \leq  c_1\exp(-c_2/\epsi^{1+2\gamma}).
    \end{equation}
  The rest of the proof is now standard.
\end{Proof}

\section{White noise approximation}
\label{sec:approximation}
In order to introduce a finite element method (FEM) that approximates
a solution of \eqref{eqn:integral.noisy.allen-cahn}, we first need to
to obtain a weak formulation in the standard sense of PDE and
FEM.  This is not possible with the presence of the white noise, so we
regularize first the problem by replacing the white noise with a
smoother stochastic term.  Our technique is inspired by that of Allen,
Novosel \& Zhang \cite{allen-novosel-zhang} for the linear heat
equation.
\subsection{A piecewise constant approximation of the white noise}
\label{sse:awn}
Consider a tensor-product partition of the space-time domain,
$\stparti$, where $\sigma,\rho\in\reals^+$ and
\begin{equation}
  \begin{gathered}
    \sparti:=\ensemble{D_m}{D_m:=(x_{m-1},x_m),\:m\in\fromto1M},
    \\
    \text{and }\tparti:=\ensemble{I_n}{I_n:=[t_{n-1},t_n),\:n\in\fromto1N},
  \end{gathered}
\end{equation}
are, respectively, a space-domain, and a time-domain, partition; each one
of these partitions is uniform, that is
\begin{equation}
  \xm-\xmo=\sigma,\:\forall m\in\fromto1M\text{ and }\tn-\tno=\rho,\:\forall
  n\in\fromto1N
\end{equation}
and $x_0=-1,\,x_M=1$, $t_0=0$ and $t_N=T$.  
We denote by $\chi_m=\one_{D_m}$ and $\fie_n=\one_{I_n}$ the
characteristic functions of the space subdomains and time subdomains
respectively.  

The (piecewise constant) \emph{approximation of white
noise}, abbreviated by {\em \awn} below, is given by the random space-time
function
\begin{equation}
  \barnoise(x,t)=\sum_{n=1}^N\sum_{m=1}^N\baretamn\chi_m(x)\fie_n(t)
\end{equation}
where the coefficients are the random variables defined by
\begin{equation}
  \label{eqn:awn-coefficients-def}
  \baretamn := \frac1{\sigma\rho}\int_I\int_D
  \chi_m(x)\fie_n(t)\d\bm(x,t).
\end{equation}
In the sequel we will use the shorthand
\begin{equation}
  \int_0^t\int_D f(x,s)\d\barbm(x,s)
  =\int_0^t\int_D f(x,s)\barnoise(x,s)\d x\d s,
\end{equation}
in spite of the integral being taken in the classical, non-stochastic, sense.
\begin{Lem}[Moments and independence of the \awn coefficients]
  \label{lem:awn-independence}
  The coefficients $\baretamn$ defined in
  \eqref{eqn:awn-coefficients-def} are \iid \gaussrv0{1/({\sigma\rho})}
  variables.
\end{Lem}
\begin{Proof}
  From the definitions of $\baretamn$ and property
  \eqref{eqn:white.noise.is.gaussian} we have
  \begin{equation}
    \begin{split}
      \baretamn
      &=\frac1{\sigma\rho}\int_I\int_D\chi_m(x)\fie_n(t)\d\bm(x,t)\\
      &=\frac1{\sigma\rho}\int_{I_n}\int_{D_m}\d\bm(x,t)
      =\frac{\bm(I_n\times D_m)}{\sigma\rho}\\
      &\in\gaussrv 0{\frac{\norm{I_n\times D_m}}{\sigma^2\rho^2}}
      =\gaussrv0{\frac1{\sigma\rho}}.
    \end{split}
  \end{equation}
To show independence compute the covariances for $m,m'\in\fromto1M$ and
$n,n'\in\fromto1N$, using \eqref{eqn:product.formula}, as follows
\begin{equation}
  \begin{split}
    (\sigma\rho)^2\EX[\baretamn\bareta_{m',n'}]
    &=\EX\left[
      \int_I\int_D\chi_m\fie_n\d\bm
      \int_I\int_D\chi_{m'}\fie_{n'}\d\bm
      \right]
    \\
    &=\int_I\int_D\chi_m(x)\chi_{m'}(x)\fie_n(t)\fie_{n'}(t)\d x\d t
    \\
    &=\delta^m_{m'}\delta^n_{n'}\sigma\rho,
  \end{split}
\end{equation}
where $\delta^i_j$ is the Kronecker symbol.
\end{Proof}

The \awn satisfies two important technical properties that we state and
prove next.

\begin{Lem}[Approximate \ito-type inequality]
  \label{lem:approximate.ito.inequality}
  For all deterministic functions $f\in\leb2(I\times D)$ the following
  holds true
  \begin{equation}
    \label{eqn:approximate.ito.inequality}
    \EX\left[\left(\int_I\int_D f(x,t)\d\barbm(x,t)\right)^2\right]
    \leq
    \int_I\int_D f(x,t)^2\d x\d t.
  \end{equation}
\end{Lem}
\begin{Proof}
  Lemma \ref{lem:awn-independence} and some manipulations yield
  \begin{equation*}
    \begin{split}
      \EX&\left[\left(\int_I\int_D f(x,t)\d\barbm(x,t)\right)^2\right]
      =
      \EX\left[
	\left(\int_I\int_D
	f(x,t)\sum_{m n}\baretamn
	\chi_m(x)\phi_n(t)\d x\d t\right)^2\right]
      \\
      &=
      \EX\left[
	\left(\sum_{m n}\baretamn\int_{I_n}\int_{D_m} 
	f(x,t)\d x\d t\right)^2\right]
      \\
      &=
      \EX\left[
	\sum_{nm}\baretamn^2\left(\int_{I_n}\int_{ D_m}
	f\right)^2
	+2\sum
	_{n\neq n',m\neq m'}
	\baretamn\bareta_{m' n}'
	\left(\int_{I_n}\int_{ D_m}f\right)
	\left(\int_{I_n'}\int_{ D_m'}f\right)
	\right]
      \\
      &=
      \sum_{m n}\EX\left[\baretamn^2\right]
      \left(
      \int_{I_n}\int_{ D_m}f
      \right)^2
      \\
      &=\sum_{m n}
      \frac1{\rho\sigma}\left(
      \int_{I_n}\int_{ D_m}f
      \right)^2
      \\
      &\leq\sum_{m n}\int_{I_n}\int_{ D_m}f^2
      =\int_I\int_D f(x,t)^2\d x\d t.
    \end{split}
  \end{equation*}
In the next-to-last step we use the Cauchy--Schwarz inequality.
\end{Proof}
\begin{Obs}
  Lemma \ref{lem:approximate.ito.inequality} and
  \eqref{eqn:ito.isometry}
  imply that
  \begin{equation}
    \label{eqn:approximate.ito.inequality.E}
    \EX\left[\left(\int_I\int_D f(x,t)\d\barbm(x,t)\right)^2\right]
    \leq
    \EX\left[\left(\int_I\int_D f(x,t)\d\bm(x,t)\right)^2\right].
  \end{equation}
  In other words, the $\leb2$-type regularity properties of the \awn
  will be, at the worse, the same as those of the white noise itself.
\end{Obs}

Since we will need bounds on space-time norms of the \awn, but in probabilities rather than in expectation, we establish the following basic result.
\begin{Lem}[$\leb\infty(\leb2)$ and $\leb2(\leb2)$ bounds for the AWN]
  \label{lem:awn-lil2-bound}
  For each $K>0$ we have
  \begin{gather}
    \label{eqn:awn-lil2-bound-prob}
    \prob\setof{\sup_{t\in[0,T]}\Norm{\barnoise(t)}_{\leb2(D)}\leq K}
    \geq
    \pospart{
    1-\frac T\rho\left(1+\frac {K^2}2\rho\right)^{1/\sigma-1}
	  \exp\left(-\frac {K^2}2\rho\right)
	  }
    \intertext{and}
    \label{eqn:awn-l2l2-bound-prob}
    \prob\setof{\Norm{\barnoise}_{\leb2(D\times[0,T])}\leq K}
      \geq
      1-\left(1+\frac{K^2}2\right)^{T/(\sigma\rho)-1}
      \exp\left(-\frac{K^2}2\right).
  \end{gather}
\end{Lem}

\begin{Proof}
  We proceed in steps.
  \begin{Steps}
  \item  Recall that $M=2/\sigma$ and $N=T/\rho$. 
    By the definition of $\barnoise$ we have, for each $t\in[0,T]$
    and $n\in\fromto1N$ such that $t\in I_n$, that
    \begin{equation}
      \label{eqn:awn-lil2-bound-basic}
	\Norm{\barnoise(t)}_{\leb2(D)}^2
	=\sigma\sum_{m=1}^M\baretamn^2
	=\frac1{\rho}\sum_{m=1}^M\etamn^2,
    \end{equation}
    where the $\etamn\in\gaussrv01$.
    In order to conclude, we will obtain a condition on the
    right-hand side that makes it smaller than $K^2$, 
    for all $n\in\fromto1N$.
  \item For each $n\in\fromto1N$ we consider
    the random variable
    \begin{equation}
      H_n := \sum_{m=1}^M\etamn^2.
    \end{equation}
    Notice that, in view of Lemma \ref{lem:awn-independence} for
    $n\neq n'$, $H_n$ and $H_{n'}$ are independent.  Let us fix $n$
    for a while and find an event for which $H_n\leq \rho K^2$.  By
    Lemma \ref{lem:awn-independence} and a basic probability fact
    \cite[Pbm. 20.16]{billingsley:book:1995}, the random variable
    $H_n$ has a chi-squared distribution with $M$ degrees of freedom.
    Its density is given by
    \begin{equation}
      \label{eqn:chi-squared-integral}
	\frac{z^{M/2-1}\exp(-z/2)}{2^{M/2}\Gamma(M/2)},
	\text{ for }z>0,
    \end{equation}
    and $0$ for $z\leq 0$, where $\Gamma$ is the Euler Gamma-function.
    Thus we have
    \begin{equation}
      \prob\setof{H_n\leq \rho K^2}
      =\frac1{2^{M/2}\Gamma(M/2)}
      \int_0^{\rho K^2}{z^{M/2-1}\exp(-z/2)}\d z.
    \end{equation}
    \item
      We prove next a lower bound on this integral in the case where
      $M$ is even, the odd case being similar.  Let $y$ play the role
      of $\rho K^2$ and consider for each $k\in\NO$ the integral
      \begin{equation}
	I_k:=\int_0^y z^k\exp(-z/2)\d z.
      \end{equation}
      An integration by parts yields the recursive expression
      \begin{equation}
	I_k = 2k I_{k-1}-2y^k\exp(-y/2),
      \end{equation}
      which allows, by an inductive argument, to see that
      \begin{equation}
	I_k = 2^{k+1}k!-2\sum_{i=0}^k\frac{k!}{(k-i)!}y^{k-i}2^i\exp(-y/2).
      \end{equation}
      An easy manipulation with the binomial formula implies that
      \begin{equation}
	I_k\geq 2^{k+1}k!\left(1-(1+y/2)^k\exp(-y/2)\right).
      \end{equation}
      Taking $k=M/2-1$ in the above and recalling the definition of
      $I_k$ and \eqref{eqn:chi-squared-integral} it follows that
      \begin{equation}
	\prob\setof{H_n\leq \rho K^2}\geq
	1-\left(1+\frac{\rho K^2}2\right)^{M/2-1}
	\exp\left(-\frac{\rho K^2}2\right);
      \end{equation}
      which implies
      \begin{equation}
	\label{eqn:awn-lil2-bound-time-loc}
	\prob\setof{H_n\leq \rho K^2}\geq
	\pospart{1-\left(1+\frac{\rho K^2}2\right)^{M/2-1}
	\exp\left(-\frac{\rho K^2}2\right)}.
      \end{equation}
    \item To conclude the proof, we introduce the event
      \begin{equation}
	\W_K^2=\bigcap_{n=1}^N\setof{H_n\leq
	\rho K^2},
      \end{equation}
      and we observe that, in view of
      \eqref{eqn:awn-lil2-bound-basic},
      on $\W_K^2$ we have
      \begin{equation}
	\Norm{\barnoise(t)}_{\leb2(D)}\leq K,
	\:\forall t\in[0,T].
      \end{equation}
      On the other hand, using the independence of $H_n$,
      $n\in\fromto1N$, the simple
      fact that $(1-\xi)^N\geq 1-N\xi$ for $\xi\leq1$ and
      \eqref{eqn:awn-lil2-bound-time-loc} we can estimate the probability
      \begin{equation}
	\begin{split}
	  \prob(\W_K^2)
	  &=\prod_{n=1}^N\prob\setof{H_n\leq \rho K^2}
	  \\
	  &\geq
	  \left(
	  \pospart{
	    1-\left(1+\frac{\rho K^2}2\right)^{M/2-1}
	    \exp\left(-\frac {\rho K^2}2\right)
	    }\right)^N
	  \\
	  &\geq
	  \pospart{
	    1-N\left(1+\frac{\rho K^2}2\right)^{M/2-1}
	    \exp\left(-\frac {\rho K^2}2\right)
	  }.
	\end{split}
      \end{equation}
      By replacing $N=T/\rho$ and $M=2/\sigma$ we get
      \eqref{eqn:awn-lil2-bound-prob}.
    \item
      Estimate \eqref{eqn:awn-l2l2-bound-prob} is obtained simply by using
      \eqref{eqn:awn-lil2-bound-time-loc} with $\rho K^2$ and $M$ replaced by
      $K^2$ and $MN$ respectively.
  \end{Steps}
\end{Proof}
\begin{Obs}[alternative proof]
  As pointed out by one of the referees, it is possible to prove Lemma
  \ref{lem:awn-lil2-bound} more directly, by using martingale inequalities.
\end{Obs}
\begin{Obs}[interpretation of 
    \eqref{eqn:awn-lil2-bound-prob} and \eqref{eqn:awn-l2l2-bound-prob}]
  We may rewrite the term appearing in
  \eqref{eqn:awn-lil2-bound-prob}, as
  \begin{equation}
    \frac T\rho\left(1+\frac {K^2}2\rho\right)^{1/\sigma-1}
	  \exp\left(-\frac {K^2}2\rho\right)
	  =:T \exp F(\rho,\sigma,K).
  \end{equation}
  A practical way to use such a result is by fixing first
  $T,\rho,\sigma\in\reals^+$ and then requiring a big enough $K$ such
  that $T \exp F(\rho,\sigma,K)\ll0$.  This is made possible by the
  fact that $\lim_{K\rightarrow\infty}F(\rho,\sigma,K)=-\infty$ for
  any fixed $\rho,\sigma\in\reals^+$.  The same type of observation is
  valid also for the \eqref{eqn:awn-l2l2-bound-prob}.
\end{Obs}

\section{The regularized solution}
\label{sec:convergence}
We now introduce the regularized solution to problem
\eqref{eqn:noisy-allen-cahn},
\eqref{eqn:initial-condition}--\eqref{eqn:boundary-conditions}, which
we obtain by replacing the white noise by the \awn in
\eqref{eqn:noisy-allen-cahn}.  The role of the regularized problem is
pivotal in devising a numerical scheme to approximate the stochastic
Allen--Cahn problem.  We discuss the approximation properties of
this regularization with respect to the original problem.
\begin{Def}[regularized solution]
The {\em regularized solution}, $\baru$, of the noisy Allen--Cahn
problem is the unique continuous solution of the integral equation
\begin{multline}
  \label{eqn:integral.representation.app}
  \baru(x,t)=-\int_0^t\int_D G_{t-s}(x,y)\fep(\baru(y,s))\d y\d s
  \\
  +\int_D G_t(x,y)u_0(y)\d y
  +\sep\int_0^t\int_D G_{t-s}(x,y)\d\barbm(y,s).
\end{multline}
\end{Def}
\begin{Lem}[maximum principle for regularized solutions]
  \label{lem:max-princ-app}
  For fixed $T,K_0>0$, there exist $\delta_0,c_1,c_2>0$, independent
  of $\epsi$, such that if $\Norm{u_0}_{\leb\infty(D)}\leq 1+\delta_0$
  then
  \begin{equation}
    \label{eqn:max-princ-app-prob}
    \prob\setof{\sup_{t\in[0,T]} \Norm{\baru(t)}_{\leb\infty(D)}\leq
      1+K_0}\geq 1-c_1\exp(-c_2/\epsi^{1+2\gamma}).
  \end{equation}
\end{Lem}
\begin{Proof}
  We follow exactly the proof of Lemma \ref{lem:prob.max.princ}, by
  observing that \eqref{eqn:approximate.ito.inequality.E} ensures that
  all the estimates for the stochastic integrals of the white noise
  can be ``translated'' in corresponding estimates for the integrals
  of the approximate white noise.  The constants appearing in this
  Theorem can be therefore taken to be the same that appear in
  \secref{lem:prob.max.princ}.
\end{Proof}
\begin{Obs}[regularized solution is strong solution]
Notice that the regularized solution $\baru$ of
\eqref{eqn:integral.representation.app} is in fact a weak solution in
the PDE sense, i.e., $\baru(t;\w)\in\sobh1(D)$ and
$\partial_t\baru(t;\w)\in\leb2(D)$ for all $t\in(0,T]$ and $\w\in\W$,
and the following \emph{weak formulation} is satisfied:
\begin{equation}
  \label{eqn:app-sac}
  \begin{split}
    \ltwop{\partial_t\baru(t;\w)}\phi+
    &\ltwop{\grad\baru(t;\w)}{\grad\phi}
    +\ltwop{\fep(\baru(t;\w))}\phi 
    \\
    &= \sep\ltwop{\barnoise(t;\w)}{\phi},
    \quad\forall\phi\in\sobhz1(D),t\in(0,T],\\
    \text{ and }\baru(0;\w)
    &=u_0,
  \end{split}
\end{equation}
for each $\w\in\W$ (the notation $\ltwop\cdot\cdot$ indicating the
inner product in $\leb2(D)$).  Indeed, each one of the \awn's
realizations, $\barnoise(\omega)$, is a piecewise constant space-time
function.  For each such realization the usual regularity theory for
semilinear parabolic equations with piecewise continuous data can be
applied and the corresponding weak formulation written down
\cite{ladyzenskaja-solonnikov-uralceva:book}.
\end{Obs}
  Our next goal is to show that the regularized approximate solution
  converges to the solution $u$.  For this we will estimate the {\em
  regularization error}
  \begin{equation}
    e(x,t)=u(x,t)-\baru(x,t),
  \end{equation}
  in terms of the white noise regularization parameters $\sigma$ and $\rho$,
  and show that it converges to zero in an appropriate sense.
\begin{The}[convergence to the stochastic solution]
  \label{the:awn-exact-convergence}
  For a fixed $T$, there exist constants $c_1$, $c_2$, $C_1$ and $C_2$
  such that for each $\epsi\in(0,1)$ there correspond an event
  $\W^\infty_\epsi$ and a constant $C_\epsi>0$ such that
    \begin{gather}
      \label{eqn:awn-exact-convergence:event}
      \prob(\W^\infty_\epsi)\geq 1 - 2 c_1\exp(-c_2/\epsi^{1+2\gamma})\text{ and }
      \\
      \label{eqn:awn-exact-convergence:rate}
      \int_{\W^\infty_\epsi}\left(\int_0^T\int_D\norm{\baru-u}^2\right)
      \d\prob\leq
      C_\epsi\left(C_1\rho^{1/2}+C_2\frac{\sigma^2}{\rho^{1/2}}\right),
      \:\forall \sigma,\rho>0.
    \end{gather}
\end{The}

\begin{Proof}
  We proceed by steps.
  \begin{Steps}
  \item By the integral representations of $u$,
    \eqref{eqn:integral.noisy.allen-cahn}, and $\baru$,
    \eqref{eqn:integral.representation.app}, we can represent the error
    as an integral too:
    \begin{equation}
      \label{eqn:regularized-noise-error}
      \begin{split}
	e(x,t)=\int_0^t\int_D G_{t-s}(x,y)
	\big(\fep(\baru(y,s))-\fep(u(y,s))\big)\d y\d s\\
	+\sep\int_0^t\int_D G_{t-s}(x,y)\left(\d\bm(y,s)-\d\barbm(y,s)\right)
      \end{split}
    \end{equation}
    for all $(x,t)\in D\times(0,T]$.  So our task now is to bound the
    terms in the left-hand side of \eqref{eqn:regularized-noise-error}
    in the appropriate norm.

  \item 
    In view of the maximum principle for both the exact solution,
    \secref{lem:prob.max.princ}, and the approximate solution,
    \secref{lem:max-princ-app}, there exists an event
    $\W^\infty_\epsi\subset\W$ such that
    \begin{gather}
      \prob(\W^\infty_\epsi)\geq 1-2c_1\exp(-c_2/\epsi^{1+2\gamma})
      \intertext{ and }
      \W^\infty_\epsi
      \subset\setof{\Norm{u(t)}_{\leb\infty(D)},
	\Norm{\baru(t)}_{\leb\infty(D)}\leq 3,\,\forall t\in[0,T]}.
    \end{gather}
    The choice of the number $3$ is quite arbitrary here.  In fact any
    number of the form $1+K_0$ with $K_0>0$ will do, with the
    appropriate change of constants.  This and the local Lipschitz
    continuity of $f$ imply that
    \begin{equation}
      \norm{\fep(\baru)-\fep(u)}
      \leq\frac{28}{\epsi^2}\norm{\baru-u},\text{ on $\W^\infty_\epsi$.}
    \end{equation}
    
  \item  Working now on the event $\W_\epsi^\infty$ and introducing
  the functions
    \begin{gather}
      \ep(r):=\int_0^r\int_D e(x,t)^2\d x\d t\\
      \phi(r):=\int_0^r\int_D\norm{\int_0^t\int_D
	G_{t-s}(x,y)(\d\bm(y,s)-\d\barbm(y,s))}^2\d x\d t
    \end{gather}
    for all $r\in[0,T]$, we infer from
    \eqref{eqn:regularized-noise-error}
    that
    \begin{equation}
      \label{eqn:pre-gronwall-1}
      \ep(r)\leq 2\int_0^r\int_D
      \left(\int_0^t\int_D\norm{G_{t-s}(x,y)}\frac{{28}}{\epsi^2}e(y,s)
      \d y\d
      s\right)^2\d x\d t
      +2\epsi^{2\gamma}\phi(r).
    \end{equation}
    The integral in \eqref{eqn:pre-gronwall-1} can be bounded,
    using the
    Cauchy--Schwarz
    inequality, by
    \begin{multline}
      2\frac{{28}^2}{\epsi^4}\int_0^r\int_D
      \left(\int_0^t\int_D\norm{G_{t-s}(x,y)}^2
      \d y\d s
      \int_0^t\int_D e(y,s)^2\d y\d s
      \right)\d x\d
      t
      =\int_0^r z(t)\ep(t)\d t
    \end{multline}
    where
    \begin{equation}
      z(t):=2\frac{{28}^2}{\epsi^4}\int_D\int_0^t\int_D\norm{G_{t-s}(x,y)}^2
      \d y\d s\d x.
    \end{equation}
    Inequality \eqref{eqn:pre-gronwall-1} implies
    \begin{equation}
      \ep(r)\leq\phi(r)+\int_0^r z(t)\ep(t)\d t,
    \end{equation}
    for each $r\in I$.  Applying
    the Gronwall lemma to this
    inequality we obtain
    \begin{equation}
      \ep(T)\leq\exp\left(\int_0^T z(t) \d t\right)\epsi^{2\gamma}\phi(T)
      \leq C_\epsi\phi(T),
    \end{equation}
    where---by estimating the heat kernel---the constant is given by 
    \begin{equation}
      C_\epsi:=\epsi^{2\gamma}\exp\left(\frac{28^2\,T}{12\,\epsi^{4}}\right).
    \end{equation}

  \item
    By summing with respect to $\prob$ on the event $\W^\infty_\epsi$ 
    both members of this inequality we obtain
    \begin{equation}
      \int_{\W^\infty_\epsi}\int_0^T\int_D
      \norm{\baru-u}^2\d x\d t\d\prob
      \leq C_\epsi\int_{\W^\infty_\epsi}\phi(T)\d\prob
      \leq C_\epsi\EX[\phi(T)].
    \end{equation}
    We conclude by observing \cite[Lem. 2.3]{allen-novosel-zhang} that
    there exist $C_1,C_2>0$, depending only on $T$, such that
    \begin{equation}
      \EX[\phi(T)]\leq C_1\rho^{1/2}+C_2\frac{\sigma^2}{\rho^{1/2}}.
    \end{equation}
    Thus we established that
    \begin{equation}
      \int_{\W^\infty_\epsi}\int_0^T\int_D\norm{\baru-u}^2
      \d x\d t\d P\leq C_\epsi\left(C_1\rho^{1/2}+C_2\frac{\sigma^2}{\rho^{1/2}}\right),
    \end{equation}
   \end{Steps}
   as we claimed.
\end{Proof}
\begin{Obs}[About the constant $C_\epsi$]
Theorem \ref{the:awn-exact-convergence} insures that, for fixed $T$
and $\epsilon$, the approximate solution $\baru$ converges to $u$ as
$\rho,\sigma\rightarrow0$.  The constant $C_\epsi$ appearing in the
estimate depends exponentially on both $1/\epsi^4$ and $T$, thus for
small $\epsi$, or large $T$, this might force us to take very small
$\rho$ and $\sigma$.  This fact should to be taken into account in
practice.  The bound we have proved seems to be pessimistic though, as
the choice of $\sigma$ and $\rho$, used in our subsequent
numerical experiments, indicates.
\end{Obs}
\begin{Obs}[convergence rate]
  \label{obs:convergence-rate}
  Recalling that $\rho$ and $\sigma$ can be thought as discretization
  (in addition to regularization) parameters, the convergence rate
  found in \eqref{eqn:awn-exact-convergence:rate} is in acccordance
  with standard results for linear parabolic equations
  \cite[e.g.]{allen-novosel-zhang}.  Note that the need for
  $\rho=C\sigma^2$, the so-called ``parabolic space-time scaling'', we
  obtain the right balance between the two terms appearing on the
  right-hand side of \eqref{eqn:awn-exact-convergence:rate}.
\end{Obs}
\section{The regularized solution's limit for $\epsi\to0$}
\label{sec:weak.noise}
In this section we focus on the relation between the regularized
stochastic Allen--Cahn problem (\ref{eqn:app-sac}) and the
deterministic version.  The reason to do this, is to find, in an
analytical setting, what conditions should be taken on the
regularization parameters, $\rho$ and $\sigma$ for the noise to be
captured in the regularized equation.

We mainly show that the the error between $\baru$ and the
deterministic solution to the Allen--Cahn problem, $q$, in an
appropriate probability-$\leb\infty(0,T;\leb2)$ sense, is of order
$\Oh(\epsi^3)$ as $\epsi\to0$ for fixed $\gamma>3$ and
$\rho,\sigma>0$.  In \secref{sse:low.noise.resolution} we give an
interpretation of this result as an evaluation of the risk of
obtaining a poor resolution of the noise for fixed $\rho$ and $\sigma$
and too small $\epsi$. This poor resolution may lead to the
disappearance of the stochastic effects in the regularized equation,
even after the appropriate rescaling, because $\baru$ becomes much
closer to $q$ than $u$, with respect to $\epsi$.  This point is
further investigated numerically in \secref{sec:numerics}.

Our proof makes use of a \emph{spectrum estimate} result for the
linearized Allen--Cahn operator
\cite{chen.xinfu:94,demottoni-schatzman:95}, which is recalled in
Theorem \ref{the:spectrum.estimate} and the $\leb2(D)$ estimate on the
noise given by \ref{lem:awn-lil2-bound}.  The proof's technique
is a continuous data dependence result for parabolic equations based
on a Bernoulli--Gronwall type argument, in the spirit of Feng \& Wu
\cite{feng-wu:05}.  

The result holds for $\gamma>3$ and it is an open
problem, as far as we know, to find the critical $\gamma$ for which
the result ceases to hold.

\subsection{Deterministic solution and resolved profiles}
Denote by $q$ be the (classical) solution of the problem
\begin{gather}
  \label{eqn:deterministic-allen-cahn}
  \partial_t q-\lap q+\fep(q)=0,\text{ in } D\times I\\
  q(0)=u_0,\text{ on }D\\
  \grad q(t,0)=\grad q(t,1)=0,\:t\in I.
\end{gather}
We consider also the function of space only $q_0$ defined as the unique solution to
\begin{equation}
  -q_0''+f_1(q_0)=0\text{ in }\reals,\quad q_0(\pm\infty)=\pm 1
  \text{ and }q_0(0)=0.
\end{equation}
(In fact it is $q_0=\tanh$.)  We will assume
from now on that $u_0$ is a \emph{resolved profile solution}, which is
defined to be an $\epsi$-linear perturbation of an $\epsi$-rescaled
and shifted $q_0$. That is, for all $x\in D$,
$u_0(x)=q_0((x-x_0)/\epsi)+\epsi p_0(x)$ where $x_0 \in D$, and $p_0$
is such that $u_0$ satisfies the Neumann boundary conditions
$u'_0(\pm1)=0$.  With this choice of initial condition the
linearization of the operator $u\mapsto-\lap u + \fep(u)$ about $q$
enjoys the following spectral property.
\begin{The}[Spectrum estimate \cite{chen.xinfu:94,demottoni-schatzman:95}]
\label{the:spectrum.estimate}
There exists a constant $\lambda_0>0$ such that for any
  $\epsi\in(0,1]$ we have
  \begin{equation}
    \label{eqn:spectrum.estimate}
    \Norm{\grad\phi}_{\leb2(D)}^2
    +\ltwop{\fep'(q)\phi}\phi
    \geq -\lambda_0\Norm{\phi}_{\leb2(D)}^2,\:\forall\phi\in\sobh1(D).
  \end{equation}
\end{The}

It is also a well-known consequence of the maximum-principle that if
$\norm{u_0}\leq1$ (which is the case when $u_0$ is a resolved profile)
then $\norm q\leq 1$.

The main result of this section is
\begin{Lem}
[continuous dependence for the regularized-deterministic error]
\label{the:weak.noise}
There exists a bounded and non-increasing function
$K_1:[0,\infty)\to\reals$ and a constant $K_2$, both depending only on
$\lambda_0$, such that
\begin{equation}
  \label{eqn:weak.noise-estimate}
  \Norm{\baru(t)-q(t)}_{\leb2(D)}\leq K_2\epsi^3
\end{equation}
provided
\begin{equation}
  \label{eqn:weak.noise-sufficient}
  \int_0^t\Norm{\barnoise(s)}_{\leb2(D)}^2\exp(-(3+2\lambda_0)s)\d s
  \leq K_1(t)\epsi^{6-2\gamma},
\end{equation}
for $t\in[0,T]$.
\end{Lem}
\begin{Proof}
  We divide the proof in several steps and we denote in it
  $\Norm{\cdot}_{\leb2(D)}$ simply by $\Norm{\cdot}$.
  \begin{Steps}
    \item We start by deriving an energy inequality for the error
      \begin{equation}
	\bare := \baru - q.
      \end{equation}
      Since $\baru$ satisfies the weak formulation \eqref{eqn:app-sac}
      and $q$ is a classical solution, we can write the following PDE
      in its weak formulation for $\bare$:
      \begin{equation}
	\label{eqn:barerror-relation}
	\ltwop{\baret}\phi+\abil{\bare}\phi
	+\ltwop{\fep'(q)\bare}{\phi}
	=\sep\ltwop{\barnoise}{\phi} 
	- \frac1{\epsi^2}\ltwop{\barr\bare^2}{\phi},
	\:\forall\phi\in\sobh1(D),
      \end{equation}
      where
      \begin{equation}
	\label{eqn:def-barr}
	\barr := 3q+\bare = 2q+\baru.
      \end{equation}
      
      Testing with $\bare$ in \eqref{eqn:barerror-relation}  
      we obtain
      \begin{equation}
	\label{eqn:basic-energy-identity}
	\ltwop{\baret}{\bare}+\Norm{\grad\bare}^2
	+\ltwop{\fep'(q)\bare}{\bare}
	\\
	\leq
	\sep\ltwop{\barnoise}{\bare}-\frac1{\epsi^2}\ltwop{\barr}{\bare^3}.
      \end{equation}      
    \item 
      The next step is to bound the terms in the right-hand side of
      \eqref{eqn:basic-energy-identity}.  The first term can be written as
      \begin{equation}
	\label{eqn:barnoise-bound}
	\sep\ltwop{\barnoise(t)}{\bare(t)}
	  \leq
	  \frac{\epsi^{2\gamma}}2\Norm{\barnoise(t)}^2
	  +
	  \frac12\Norm{\bare(t)}^2.
      \end{equation}

      To produce a bound on the second term of the right-hand side of
      \eqref{eqn:basic-energy-identity} we use \eqref{eqn:def-barr},
      valid in $1$ spatial dimension, to obtain
      \begin{equation}
	\ltwop{\barr(t)}{\bare(t)^3}
	=3\ltwop q{\bare(t)^3}+\Norm{\bare(t)}_{\leb4(D)}^4.
      \end{equation}
      By the fact that $\norm q\leq 1$ and the Sobolev embedding
      $\sobh1(D)\embeddedin\leb\infty(D)$, valid for $D\subset\reals$,
      the first term on the right hand side can be bounded using
      \begin{equation}
	\begin{split}
	  \norm{3\ltwop{q}{\bare(t)^3}}
	  &\leq
	  3\Norm{\bare(t)}_{\leb\infty(D)}\Norm{\bare(t)}^2
	  \leq 
	  C_1\Norm{\bare(t)}_{\sobh1(D)}\Norm{\bare(t)}^2
	  \\
	  &\leq 
	  {\lambda_1\epsi^4}\Norm{\bare(t)}_{\sobh1(D)}^2
	  +\frac{C_1^2}{4\lambda_1\epsi^4}\Norm{\bare(t)}^4
	\end{split}
      \end{equation}
      where $C_1$ is $3$ times the Sobolev embedding constant for $D$
      and $\lambda_1:=\min\setof{1,\lambda_0}$ (the reason for this choice
      will be apparent in the next step) with $\lambda_0$ from
      \eqref{eqn:spectrum.estimate}.  As a consequence we have
      \begin{equation}
	\label{eqn:bare-nonlin-bound}
	-\frac1{\epsi^2}\ltwop\barr{\bare^3}
	\leq
	\lambda_1\epsi^2\Norm{\bare}_{\sobh1(D)}^2
	+\frac{C_2}{\epsi^6}\Norm{\bare}^4
	-\frac1{\epsi^2}\Norm{\bare}_{\leb4(D)}^4
      \end{equation}
      where $C_2=C_1^2/(4\lambda_1)$.
      
    \item
      Owing to the spectrum estimate \eqref{eqn:spectrum.estimate} and
      the fact that $f'(q)\geq-1$ we have
      \begin{equation}
	\label{eqn:energy-spectrum-below}
	\begin{split}
	  \Norm{\grad\bare}^2+\ltwop{\fep'(q)}{\bare^2}
	  &=:A=(1-\epsi^2)A+\epsi^2A
	  \\
	  &\geq-(1-\epsi^2)\lambda_0\Norm{\bare}^2
	  +\epsi^2\Norm{\grad\epsi}^2-\Norm{\bare}^2
	  \\
	  &=-\left((1-\epsi^2)\lambda_0+1\right)\Norm{\bare}^2
	  +\epsi^2\Norm{\grad\bare}^2
	  \\
	  &=-(1+\lambda_0)\Norm{\bare}^2
	  +\epsi^2(\Norm{\grad\bare}^2+\lambda_0\Norm{\bare}^2)
	  \\
	  &\geq-(1+\lambda_0)\Norm{\bare}^2
	  +\lambda_1\epsi^2\Norm{\bare}_{\sobh1(D)}^2.
	\end{split}
      \end{equation}

      The inequalities (\ref{eqn:basic-energy-identity}),
      (\ref{eqn:barnoise-bound}), (\ref{eqn:energy-spectrum-below})
      and (\ref{eqn:bare-nonlin-bound}) lead to
      \begin{multline}
	\label{eqn:pregronwall-a}
	  \frac12\d_t\Norm{\bare(t)}^2-(1+\lambda_0)\Norm{\bare(t)}^2
	  \\
	  \leq
	  \frac{\epsi^{2\gamma}}2\Norm{\barnoise(t)}^2
	  +\frac12\Norm{\bare(t)}^2
	  +\frac{C_2}{\epsi^6}\Norm{\bare(t)}^4
	  -\frac1{\epsi^2}\Norm{\bare(t)}_{\leb4(D)}^4
	  \:\forall t\in I.
      \end{multline}

      Consider, for the rest of the proof, the following notation:
      \begin{equation}
	g(t):=\Norm{\bare(t)}^2,
	\qquad 
	a:=(3+2\lambda_0),
	\qquad
	b:=2C_2/\epsi^6,
	\qquad
	r(t):=\epsi^{2\gamma}\Norm{\barnoise(t)}^2-\Norm{\bare(t)}_{\leb4(D)}^4.
      \end{equation}
      Then \eqref{eqn:pregronwall-a} implies
      \begin{equation}
	\label{eqn:prebernoulli-1}
	g'(t)\leq a g(t)+b g(t)^2+r,\:\forall t\in[0,T].
      \end{equation}

    \item
      To proceed we will apply a Bernoulli differential inequality
      technique, which generalizes the Gronwall Lemma, in order to get
      a bound on $g(t)$.  We follow Feng \& Wu \cite[Lem.2.1]{feng-wu:05}

      Fix a $t\in[0,T]$ and let
      \begin{equation}
	\ro(s):=\int_0^s\exp(-a\tau)r(\tau)\d\tau
	\text{ and }
	p(s):=p_t(s)
	=(\ro(t)-\ro(s))\exp(as),
      \end{equation}
      for all $s\in[0,t]$.  Since
      \begin{equation}
	p'(s)=\d_s p_t(s)=-r(s)+a p(s)\text{ and }p(s)\geq0,\text{ for }s\in[0,t],
      \end{equation}
      we may write
      \begin{equation}
	\d_s(g(s)+p(s))\leq
	a(g(s)+p(s))+b(g(s)+p(s))^2,\text{ for }s\in[0,t].
      \end{equation}
      Introducing $z(s):=1/(g(s)+p(s))$,
      we can rewrite this inequality as
      \begin{equation}
	z'(s)+a z(s)\geq-b,\:\forall s\in[0,t].
      \end{equation}
      Multiplying by $\exp(as)$ and integrating over $[0,t]$ we obtain
      \begin{equation}
	z(t)\geq z(0)\exp(-at)-\frac{b(1-\exp(-at))}a.
      \end{equation}
      Noting that $g(0)=\Norm{\bare(0)}=0$, $p_t(0)=\ro(t)$ and
      $p_t(t)=0$, this yields
      \begin{equation}
	\frac1{g(t)}\geq\frac{a-b\ro(t)(\exp(at)-1)}{a\exp(at)\ro(t)}.
      \end{equation}
      We now invert both sides of this inequality, under the
      sufficient condition that
      \begin{equation}
	\label{eqn:bernoulli-sufficient}
	a-b\ro(t)(\exp(at)-1)\geq0,
      \end{equation}
      and we get
      \begin{equation}
	\Norm{\bare(t)}^2\leq\frac{a\exp(at)\ro(t)}{a-b\ro(t)(\exp(at)-1)}.
      \end{equation}
    \item
      To conclude we want to interpret more explicitly this
      result. Let us replace first \eqref{eqn:bernoulli-sufficient} by
      the sufficient condition
      \begin{equation}
	a-b\ro(t)(\exp(at)-1)\geq\delta(t),
      \end{equation}
      for some $\delta(t)>0$ that will be chosen appropriately.  This
      is equivalent to
      \begin{equation}
	\label{eqn:bernoulli-sufficient-delta}
	\int_0^t\exp(-as)r(s)\d s\:(=\ro(t))\:
	\leq\frac{a-\delta(t)}{b(\exp(at)-1)}.
      \end{equation}
      This can be ensured if we assume
      \begin{equation}
	\int_0^t\Norm{\barnoise(s)}^2\exp(-as)\d s
	\leq
	\frac{(a-\delta(t))\epsi^{6-2\gamma}}{2C_2(\exp(at)-1)}.
      \end{equation}

      Under this condition we obtain the bound
      \begin{equation}
	\label{eqn:postbernoulli}
	\Norm{\bare(t)}^2\leq
	\frac{a(a-\delta(t))\exp(at)}{2C_2\delta(t)(\exp(at)-1)}\epsi^6.
      \end{equation}

    \item
      We conclude by taking
      \begin{equation}
	\delta(t):=\max\setof{a-\exp(at)+1,\frac a2},
      \end{equation}
      i.e.,
      \begin{equation}
	\delta(t):=\begin{cases}
	a-\exp(at)+1,&\text{ for }t\leq t_a
	\\
	a/2,&\text{ for }t>t_a
	\end{cases}
      \end{equation}
      where $t_a=\log(1+a/2)/a$.
      Then, after putting
      \begin{equation}
	K_1(t):=\frac{\min\setof{1,a/2(\exp(at)-1)}}{2C_2}\text{ for }t\geq0,
      \end{equation}
      condition \eqref{eqn:bernoulli-sufficient-delta} may be replaced by
      \begin{equation}
	\int_0^t\Norm{\barnoise(s)}^2\exp(-as)\d s
	\leq
	K_1(t)\epsi^{6-2\gamma},
      \end{equation}
      a condition under which we have, from \eqref{eqn:postbernoulli}
      \begin{equation}
	\Norm{\bare(t)}^2\leq
	{K_2}^2\epsi^6,
      \end{equation}
      where ${K_2}^2:=\sup_{[0,\infty)}K_1(t)^2\exp(at)<\infty$.
  \end{Steps}
\end{Proof}

As a consequence of this estimate we state the following result,
which, roughly speaking, implies that in order for the noise to have the
chance of persisting in the limit, as $\epsi\to0$, the parameters
$\rho,\sigma$ must also go to zero.
\begin{The}[low-intensity approximate white noise]
\label{the:low.noise}
There exists a constant $C=C(\lambda_0)$ such that for all fixed
$\gamma>3,\rho,\sigma,T>0$ we have
\begin{equation}
  \label{eqn:low.noise}
  \lim_{\epsi\to0}\prob\ensemble{\w\in\W}
      {\sup_{[0,T]}\Norm{\baru(\cdot;\w)-q}_{\leb2(D)}<C\epsi^3}
      =1.
\end{equation}
\end{The}
\begin{Proof}
  Let $C=K_2$ in \eqref{the:weak.noise}.  Choosing to use an
  $\leb2(0,t;\leb2(D))$-norm estimate for the white noise, for the
  estimate \eqref{eqn:weak.noise-estimate}, in view of the
  monotonicity of $K_1$, it is enough to assume the sufficient
  condition
  \begin{equation}
    \label{eqn:bernoulli-sufficient-l2l2}
    \int_0^t\Norm{\barnoise(s)}^2\d s
    \leq K_1(T)\epsi^{6-2\gamma}.
  \end{equation}
  According to \eqref{eqn:awn-l2l2-bound-prob}, this condition is
  satisfied with probability
  \begin{equation}
    1-\left(1+\frac{K_1(T)}2\epsi^{6-2\gamma}\right)^{T/(\sigma\rho)-1}
    \exp\left(-\frac{K_1(T)}2\epsi^{6-2\gamma}\right).
  \end{equation}
  Since $\gamma>3$, $6-2\gamma<0$ and, for fixed $\sigma,\rho$ and
  $T>0$, it is possible to make this probability arbitrarily close to
  $1$ for $\epsi>0$ small enough, as claimed.
\end{Proof}

\begin{Obs}
  Note that, it is also possible to obtain a variant of
  \ref{the:low.noise} by employing \eqref{eqn:awn-lil2-bound-prob}
  instead of \eqref{eqn:awn-l2l2-bound-prob}.  This leads to slightly
  better lower estimates of the probability for the same $\epsi>0$ for
  longer time $T$.
\end{Obs}

\subsection{White noise resolution by the \awn}
  \label{sse:low.noise.resolution} 
  We now describe one interpretation Theorem \ref{the:low.noise}.

  Note firstly that, in the theorem's statement, the \awn regularization
  parameters $\sigma$ and $\rho$ are kept fixed while $\epsi\to0$.  It
  is well known that, if $u_0$ is a resolved profile with center at
  $x_0$ (see \secref{sse:centers} for a definition of center), i.e.,
  $u_0=\tanh((x-x_0)/2\epsi)$, as $\epsi\to0$, the solution
  $q=q^\epsi$ of \eqref{eqn:deterministic-allen-cahn} converges, in an
  appropriate sense, to the stationary \emph{step function}
  $\chi_{x_0}:=\one_{(x_0,\infty)}-\one_{(-\infty,x_0)}$.  So Theorem
  \ref{the:low.noise} is saying that for fixed $\rho,\sigma>0$ and for
  $\epsi\to0$, the solution, $\baru$, of the regularized problem
  \eqref{eqn:app-sac} converges to this stationary step function.

  Second we note that, owing to a result by Funaki \cite[Theorem
  8.1]{funaki:95} or a similar one by Brassesco et
  al. \cite{brassesco-demasi-presutti:95}, there exists a stochastic
  process $(t,\w)\mapsto\xi^\epsi_t(\w)$ such that
  \begin{equation}\label{eqn:funaki.convergence}
    \lim_{\epsi\searrow0}\prob
    \ensemble{\w\in\W}{\sup_{t\in[0,T\epsi^{-1-2\gamma}]}\Norm{u(\cdot,t;\w) -
	\chi_{\xi^\epsi_t(\w)}}_{\leb2(D)}>\delta}=0,
  \end{equation}
  for each fixed $\delta>0$, where $\chi_{x_0}$ is the step function defined
  above, and $(t,\w)\mapsto\xi^\epsi_t$ is a stochastic process which
  converges as $\epsi\to0$, in an appropriate sense (in law), to the
  standard Brownian motion rescaled as to have diffusion coefficient
  $\sqrt{c_0}\epsi^{1/2+\gamma}$ where $c_0=3\sqrt2/4$. Of
  course, as $\epsi\to0$, $\xi^\epsi_t(\w)\to x_0$, where $x_0$ is the
  center of the initial condition $u_0$; this implies that the limit
  of $u$ and $\baru$ are consistent, as $\epsi\to0$, even when $\rho$
  and $\sigma$ are kept fixed.

  Suppose now that one wishes, in view of Theorem
  \ref{the:awn-exact-convergence}, to use the regularized solution
  $\baru(t)$, instead of $u(t)$, to approximate the diffusion
  coefficient of the process $t\mapsto{\xi^\epsi_t}$.  One way of
  doing this would be to approximate (numerically) $\baru(t,\w)$, for
  $\w\in\W$ (or a discrete analog), find its center, if it exists,
  $\bar\xi^\epsi_t(\w)$, and finally compute its average (excluding
  solutions that have no center) and its variance over $\w\in\W$. The
  resulting variance, rescaled appropriately, i.e.,
  $t\mapsto\var[\xi^\epsi_t]/(\epsi^{1+2\gamma})$, for the computation
  to be meaningful, one should see, asymptotically as $\epsi\to0$, a
  linear function $t\mapsto c_0t$.  This rescaling, which we shall
  call the Mueller--Funaki rescaling \cite{funaki:95}, is necessary in
  order to get a result that is essentially independent of $\epsi$ and
  thus easy to visualize.

  Theorem \ref{the:low.noise} tells us that, for a fixed
  $\rho,\sigma>0$, the rate of convergence of $\baru\to q$ is
  $\Oh(\epsi^3/2)$.  Since the distance
  $\left(\EX\Norm{q(t)-\chi_{\xi^\epsi_t}}^2\right)^{1/2}$, as can be
  seen using a piecewise constant approximation of $\tanh$, is
  $\Oh(\epsi^{1/2})$, it follows that $\baru$ is closer to $q$ than
  $\chi_{\xi^\epsi_t}$ and that any statistics conducted on $\baru$
  may lead to wrong results.  This is a strong indication, which is
  confirmed by the numerical results in \secref{sec:numerics}, that in
  order to capture the stochastic effects the parameters $\sigma$ and
  $\rho$ must be chosen as functions of $\epsi$.  Note that a similar
  conclusion can be derived from Theorem
  \ref{the:awn-exact-convergence} in case the dependence of $C_\epsi$
  proves to be effective, but the nature of this similar conclusion
  has its roots in deterministic considerations rather than stochastic
  ones.

  Although we have proved Theorem \ref{the:low.noise} for values of
  $\gamma>3$, it is natural to expect similar results for lower values
  of $\gamma$. In fact, our numerical experiments in
  \ref{sec:numerics} indicate that this is the case.

\section{An Euler-Galerkin finite element scheme}
\label{sec:fem}
We introduce now the finite element discretization of the regularized problem 
\eqref{eqn:app-sac}.
\subsection{Discretization partitions}
We begin by introducing the space and time partitions
\begin{equation}
  \begin{gathered}
    \shparti:=\ensemble{D'_m}{D'_m:=(x'_{m-1},x'_m),\:m\in\fromto1{M'}},
    \\
    \text{and }\thparti:=\ensemble{I'_n}{I'_n:=[t'_{n-1},t'_n),\:n\in\fromto1{N'}}.
  \end{gathered}.
\end{equation}
These partitions do not necessarily coincide with the partitions
$\sparti$ and $\tparti$ used for the regularization procedure in
\secref{sse:awn}.  Bearing in mind that this setting could be further
generalized, we limit ourselves here to the case where the numerical
discretization partitions, $\shparti$ and $\thparti$, are refinements
of the white noise regularization partitions $\sparti$ and $\tparti$,
respectively.  For each $D'_m\in\shparti$ there exists
$D_l\in\sparti$ such that $D'_m\subset D_l$ etc; this determines a
unique mapping $\mu:\fromto0{M'}\rightarrow\fromto0M$, such that
$D'_m\subset D_{\mu(m)}$.  For simplicity, we also assume that the
partitions are uniform and that the \emph{meshsize} and
\emph{timestep} are denoted respectively by $h$ and $k$.  The reason
we do not make these partitions coincide is that for the finite
element method's convergence analysis it may prove useful to have
more involved couplings of the type $h=h(\sigma)$ and $k=k(\rho)$.
In this article we consider only the simplest situation
possible where $h=\sigma$ and $k=\rho$.

\subsection{Finite element space and the discrete scheme}
Let $\fes{}\subset\honezd$ be the space of continuous piecewise linear
functions associated with the partition $\shparti$, we define the
\emph{(spatial) semi-discrete solution} as the time-dependent random finite
element function $U:[0,T]\times\W\rightarrow\fes{}$ 
which solves the SDE
\begin{equation}
  \label{eqn:fem-noisy-allen-cahn}
  \ltwop{\partial_t U(t)} V+\abil {U(t)}V+\ltwop{\fep(U(t))}V=\ltwop\barnoise V,
  \:\forall V\in\fes{},\,t\in[0,T].
\end{equation}
We discretize further this SDE in the time variable by taking a
semi-implicit Euler scheme in time associated to the partition
$I=\{t_0\}\cup\bigcup_m I_m$
\begin{equation}
  \label{eqn:semiimplicit-fem-sac}
    \ltwop{\frac{U^n-U^{n-1}}{k}}V+\abil{U^n}V
    +\ltwop{\fep(U^{n-1})}V =\ltwop\barnoise V,
    \:\forall V\in\fes{},\,t\in[0,T].
\end{equation}
The scheme is called semi-implicit, in that the linear part is treated
implicitly while it is explicit in the nonlinearity.  This means that
at each timestep only one linear problem has to be solved and no nonlinear
solver is needed.

In practice, it is more practical to use a modified version
of \eqref{eqn:semiimplicit-fem-sac} given by
\begin{equation}
  \label{eqn:newton-fem-sac}
  \begin{split}
    \ltwop{\frac{U^n-U^{n-1}}{k}}V+\abil{U^n}V
    +\ltwop{\fep'(U^{n-1})U^n}V\\
    =\ltwop{\fep'(U^{n-1})U^{n-1}-\fep(U^{n-1})}V 
    +\sep\ltwop\barnoise V,
    \:\forall V\in\fes{},\,t\in[0,T].
  \end{split}
\end{equation}
which allows bigger timesteps $k$
\cite{kessler-nochetto-schmidt:03}.  Note that this amounts to a
linearization involving one step of the Newton method to solve the
nonlinear (fully implicit) backward Euler scheme.

\subsection{The linear time-stepping system}
Let us indicate the basis functions of $\fes{}$ by $\Phi_m$, for
$m\in\fromto0{M'}$; that is the piecewise linear continuous function
such that $\Phi_m(x_l)=\delta^m_l$, for $l\in\fromto0{M'}$.  If we
indicate by $\vec u^n=(u^n_m)$ the vector of nodal values
corresponding to the discrete solution $U^n$ at time $t^n$, that is
$U^n(x)=\sum_{m=0}^{M'}u^n_m\Phi_m(x)$, then, we can translate
\eqref{eqn:newton-fem-sac} in the following matrix form
\begin{equation}
  \label{eqn:newton-fem-sac-matrix}
  \left[\frac1k\vec M+\vec A
    +\frac1{\epsi^2}\vec N(\vec u^{n-1})\right]\vec u^n
  =
  \frac1{\epsi^2}\vec g(\vec u^{n-1}) + \frac1k\vec M\vec u^{n-1} 
  + \sep\vec w,
\end{equation}
where $\vec M$, $\vec A$ are the usual finite element mass and
stiffness matrices, respectively, $\vec N(\vec u^{n-1})$ and $\vec
g(\vec u^{n-1})$ are a ``nonlinear'' mass matrix and load vector,
respectively, and $\vec w=(w_m)$ a random load vector generated at
each time-step.  A short calculation shows that for $m$ an internal
degree of freedom (node) we have
\begin{equation}
\label{eqn:random-vector-comp}
w_m = \frac{h}{2\sqrt{\sigma\rho}}(\eta_{\mu(m)-1}+\eta_{\mu(m)})
\end{equation}
where $\mu$ is the mapping introduced earlier in this section and
$\eta_l$ is a $\gaussrv01$ random number for $l\in\fromto0M$, or zero
for $l=-1,M+1$ (the boundary cases).  If the partitions $\sparti$ and
$\shparti$ coincide, which will be the case in the next section, then
$h=\sigma$ and \eqref{eqn:random-vector-comp} simplifies to
\begin{equation}
  \label{eqn:random-vector-comp.balanced}
  w_m = \frac12\sqrt{\frac{h}{\rho}}(\eta_{m-1}+\eta_m).
\end{equation}
This is the form that we employ in our computations below.  From now
on, we will consider $\tparti$ and $\thparti$ to coincide, i.e.,
$k=\rho$ in all our computations.
\section{Computations} 
\label{sec:numerics}
\newcommand{\picwidth}{.46\textwidth}
\newcommand{\SCALE}{0.75}
\newcommand{\GA}{0}
\newcommand{\GAM}{0\GA}
\newcommand{\PGAM}{0.\GA}
\newcommand{\captionA}{
  Interface average (top), $\EX[\Xi^n]$, 
  and its rescaled variance
  (bottom), $\var[\Xi^n]/20\epsi^{1+2\gamma}$, 
  as functions of
  time $t_n$, for $\epsi=0.04$.}
\newcommand{\captionB}{
  Interface average (top), $\EX[\Xi^n]$, and its rescaled
  variance (bottom), $\var[\Xi^n]/20\epsi^{1+2\gamma}$, as
  functions of time $t_n$, for $\epsi=0.02$.
}
\newcommand{\captionC}{
  Interface average (top), $\EX[\Xi^n]$, and its rescaled variance
  (bottom), $\var[\Xi^n]/20\epsi^{1+2\gamma}$, as functions of
  time $t_n$, for $\epsi=0.01$.
}
\newcommand{\captionD}{
  Each of these graphs show $\log\var[\Xi^n]$ versus $\log\epsi$ for
  the finest refinement level $l=9$. Each graph corresponds to a fixed
  choice of the discrete time index, $n$, among a sequence
  ${n_1}<{n_2}<\ldots<{n_I}$, with $n_i=10 n_{i-1}$.  A reference
  line of slope $1+2\gamma$ which is the slope that is expected for
  the other lines, the code to pick up.
}
\begin{figure}[t]
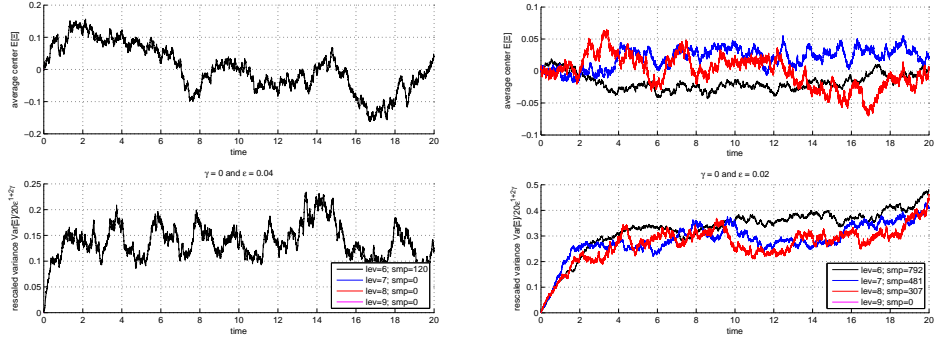
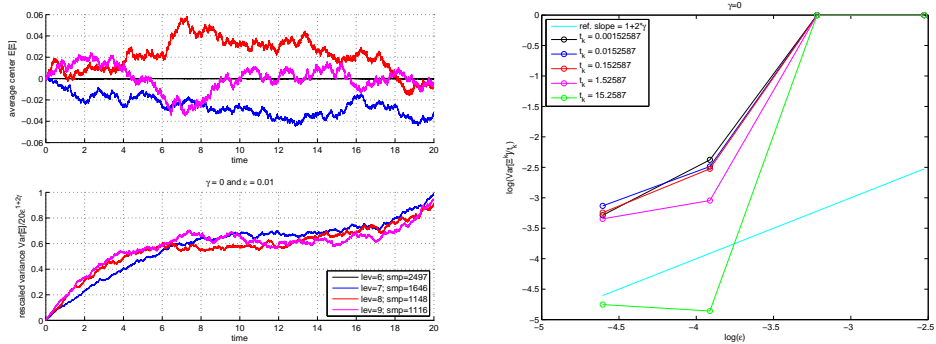

  \begin{center}
    \subfigure[{
	\label{fig:gamma\GAM-epsilon004}
	\captionA\ 
	Only one line, corresponding to the coarsest refinement level
	$l=6$, i.e., meshsize $h=2^{-5}$, is visible as for all the
	other levels no admissible paths survive up to time $T=20$
	with a center (i.e., only one zero).
    }]{
      \includegraphics[width=\picwidth,trim= 0 0 0 0,clip,scale=\SCALE]%
		      {Pictures/sac-int-004-\GAM.\graphext}
    }
    \phantom{m}
    \subfigure[{
	\label{fig:gamma\GAM-epsilon002}
	\captionB\
	Comparing this with the case $\epsi=0.004$ (Figure
	\ref{fig:gamma\GAM-epsilon004}) we see that more graphs are
	visible in each diagram, corresponding to the refinement
	levels $l=6,\,7,\,8$. This is due the fact that for smaller
	value of $\epsilon$ more sample paths the noise affects the solution 
	less.
      }
    ]{
      \includegraphics[width=\picwidth,trim= 0 0 0 0,clip,scale=\SCALE]%
		      {Pictures/sac-int-002-\GAM.\graphext}
    }
    \\
    \subfigure[{
	\label{fig:gamma\GAM-epsilon001}
	\captionC\
	In each plot here, since $\epsi$ is even smaller, we see four
	graphs corresponding to all four refinement levels tested,
	$l=6,\ldots,9$.  Note that the variance for $l=6$ is, practically
	zero, which is far from the exact results.  We believe,
	that this is due to two concurrent factors: (1) the meshsize is too
	coarse to resolve the interface layer, (2) the meshsize is too
	coarse to resolve properly the noise and the numerical
	solution is closer to the exact deterministic solution rather
	than the stochastic one.  Note also that the variance time-dependence, which is
	initially linear as expected, starts degenerating at about $T=4$.
    }]{
      \includegraphics[width=\picwidth,trim= 0 0 0 0,clip,scale=\SCALE]%
		    {Pictures/sac-int-001-\GAM.\graphext}
    }
    \phantom{m}
    \subfigure[{
	\label{fig:gamma\GAM-loglog}
	These are the results from a different perspective, with a
	two-fold purpose: (1) to check the scaling power $1+2\gamma$
	in the diffusion coefficient, and (2) to see how the quality
	of results depends behaves with respect to time. \captionD\  We see that this is reflected by the
	computations, but with the quality worsening as time grows.
    }]{
      \includegraphics[width=\picwidth,trim= 0 0 0 0,clip,scale=\SCALE]%
		      {Pictures/sac-int-\GAM-loglog.\graphext}
    }
  \end{center}
  \caption
      [Numerical results for $\gamma=\PGAM$] 
      {Numerical results for $\gamma=\PGAM$ \label{fig:gamma\GAM}}
\end{figure}
\renewcommand{\GA}{2}
\renewcommand{\GAM}{0\GA}
\renewcommand{\PGAM}{0.\GA}
\begin{figure}[t]
  \begin{center}
    \subfigure[{
	\label{fig:gamma\GAM-epsilon004}
	\captionA\ 
	Compared to Figure \ref{fig:gamma00-epsilon004} we see here
	that, due to the more modest intensity of the noise, three
	refinement levels produce enough admissible sample paths.
    }]{
      \includegraphics[width=\picwidth,trim= 0 0 0 0,clip,scale=\SCALE]%
		    {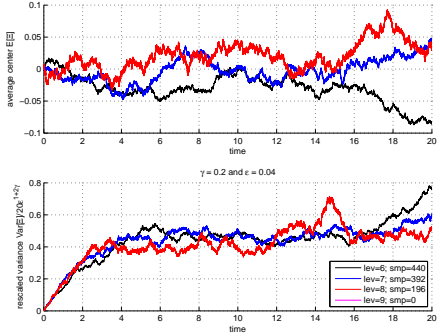}
    }
    \phantom{m}
    \subfigure[{
	\label{fig:gamma\GAM-epsilon002}
	\captionB\
	All four refinement levels are represented here.  A clear
	linear trend in the computed variance, up to time $10$ (with a
	minor deterioration between $10$ and $20$) is visible.  This
	linear dependence is the expected behavior from the theory
	saying that $\Var[\Xi^n]\approx c_0\epsi^{1+2\gamma}t_n$.
    }]{
      \includegraphics[width=\picwidth,trim= 0 0 0 0,clip,scale=\SCALE]%
		    {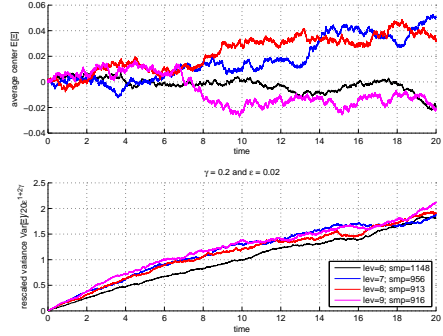}
    }
    \\
    \subfigure[{
	\label{fig:gamma\GAM-epsilon001}
	\captionC\
	We see that the linear dependence on time is clearly visible
	in the variance, but for low refinement levels (coarse
	meshsize) the slope may not be the proper one.  Also for very
	low refinement, the stochastic dynamics are not picked up at
	all.
    }]{
      \includegraphics[width=\picwidth,trim= 0 0 0 0,clip,scale=\SCALE]%
		    {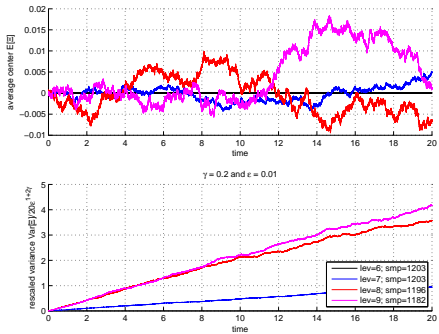}
    }
    \phantom{m}
    \subfigure[{
	\label{fig:gamma\GAM-loglog}
	\captionD\
	Comparing with the case $\epsi=0.0$, an improvement in the
	quality of the behavior as $t_n$ grows is visible, but still a
	clear deterioration in time is present.
    }]{
      \includegraphics[width=\picwidth,trim= 0 0 0 0,clip,scale=\SCALE]%
		    {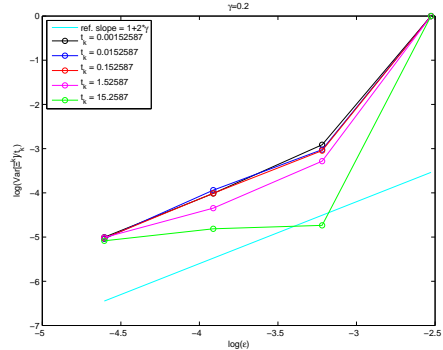}
    }
  \end{center}
  \caption
      [Numerical results for $\gamma=\PGAM$] 
      {
	\label{fig:gamma\GAM}
	Numerical results for $\gamma=\PGAM$.  A comparison with
	Figure \ref{fig:gamma00} is due and leads to two important
	observations. (1) The statistics are more robust in time
	(especially for $\epsi$ small) in that the variance shows a
	clear linear dependence on $t$ in Figures
	\ref{fig:gamma\GAM-epsilon002} and
	\ref{fig:gamma\GAM-epsilon001} and the correct scaling
	behavior in Figure \ref{fig:gamma\GAM-loglog}. This is due to
	the lower intensity of the noise because $\gamma$ is
	bigger. (2) On the other hand, we must work harder with the
	refinement level in order to pick up the stochastic effects.
	We believe that this is a practical aspect of the discussion
	in \secref{sse:low.noise.resolution}.  This poor performance for
	bigger meshsizes could not directly related to the interface
	layer resolution (a deterministic effect) as for the same
	values of $\epsi$ and $h$ but lower values of $\gamma$ we
	obtain meaningful statistics, as shown in Figure
	\ref{fig:gamma00}.
      }
\end{figure}
\renewcommand{\GA}{5}
\renewcommand{\GAM}{0\GA}
\renewcommand{\PGAM}{0.\GA}
\begin{figure}[t]
  \begin{center}
    \subfigure[{\captionA\ }]{
      \includegraphics[width=\picwidth,trim= 0 0 0 0,clip,scale=\SCALE]%
		      {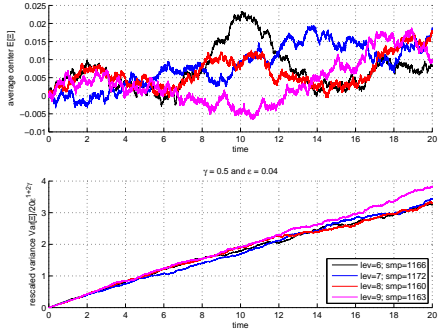}
    }
    \phantom{m}
    \subfigure[{\captionB\ }]{
      \includegraphics[width=\picwidth,trim= 0 0 0 0,clip,scale=\SCALE]%
		      {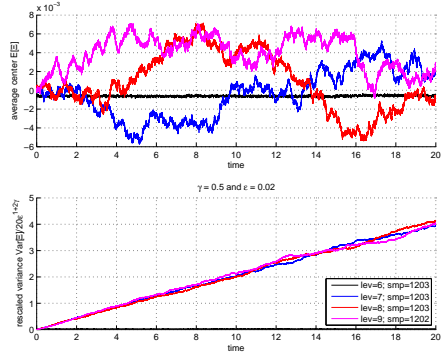}
    }
    \\
    \subfigure[{
	\captionC\ 
	The linear time dependence of the (rescaled) variance is clear
	here.  Note that for level of mesh refinement $l=6,7$ the
	computed solution is basically deterministic (average is $0$
	and variance is $0$).  This is a further indication of the
	practical importance of noise resolution, following the observations in
	\secref{sse:low.noise.resolution}.
    }]{
      \includegraphics[width=\picwidth,trim= 0 0 0 0,clip,scale=\SCALE]%
		    {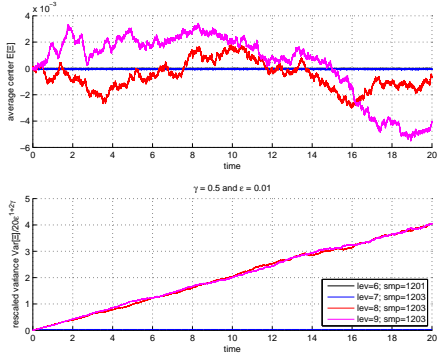}
    }
    \phantom{m}
    \subfigure[{\captionD\ 
	The value of $\gamma$ is quite high now as to keep the
	statistics robust with respect to time.
    }]{
      \includegraphics[width=\picwidth,trim= 0 0 0 0,clip,scale=\SCALE]%
		      {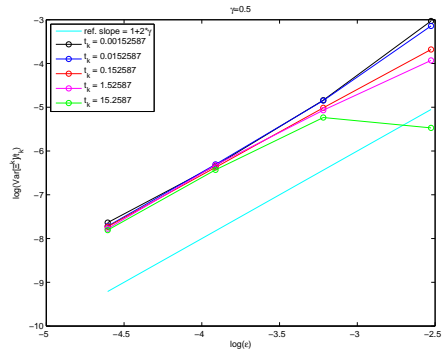}		    
    }
  \end{center}
  \caption
      [Numerical results for $\gamma=\PGAM$] 
      {
	\label{fig:gamma\GAM} 
	Numerical results for $\gamma=\PGAM$.  The same comments made
	in Figure \ref{fig:gamma02} apply here, but with an even
	clearer linear dependence of the variance upon time and a
	clearer need for high mesh refinement in order to resolve the
	white noise.  Also the computations are more robust with
	respect to lower values of $\gamma$ and they show almost no
	deterioration in time.
      }
\end{figure}

We discuss now numerical simulations of the Allen--Cahn problem, using
the scheme \eqref{eqn:newton-fem-sac-matrix}.  Our main purpose is to
match the behavior of the exact stochastic solution with the behavior
of the Monte-Carlo type numerical solution.  By ``exact behavior'' we
refer to the theoretical results of Funaki \cite{funaki:95} and
Brassesco et al \cite{brassesco-demasi-presutti:95}.  Our numerical
experiments are directed towards testing relations between the various
parameters $\epsi,\gamma$ and $\rho,\sigma$.  In addition to taking
$k=\rho$, we also take $h=\sigma$.

\subsection{Monte Carlo simulations}
We use the finite element scheme \eqref{eqn:newton-fem-sac} in
combination with Monte Carlo type simulations.  For each choice of
parameters $\epsi$, $\gamma$ and meshsize $h$, we choose the timestep
$k=h^2$ and compute between $1000$ and $2500$ samples paths, each with
a different seed for the pseudo-random number generator.

Each sample path runs from time $0$ to a final time $T$.  At the
beginning of each run, the random number generator is seeded and the
subsequent $\eta_m$ appearing in
\eqref{eqn:random-vector-comp.balanced} are chosen according to this
seed for all the run.  The seed for each sample path is determined by
the clock of the machine at the start of each run (these are also
recorded for rerunning purposes).  In this section we denote the
numerical solution (which tacitly depends on $\epsi$, $\gamma$, $h$,
etc.) by $(U^n_\w)_{n\in\fromto0N}$, where $n$ corresponds to the
timestep and $\w$ is a discrete sample, i.e., the choice of the
initial seed. We indicate by $\bar\W$ the discrete sample space,
which can be thought of being the set of all initial seeds employed.

\subsection{Phase-separation interfaces and centers}
\label{sse:centers}
Our benchmarking procedure consists in comparing the behavior of the
center of the discrete solution with that of the exact solution.  A
function $v\in\sobh1(D)$ is said to have a center if $v(x)=0$ is
uniquely solvable with a sign change for $x\in D=(-1,1)$; the solution
of $v(x)=0$ is called \emph{center} of $v$.  For example the function
$u_0(x)=\tanh((x-x_0)/\sqrt2\epsi)$ has center $x_0$.

The center represents a $0$-dimensional interface separating two
phases in the $1$-dimensional Allen--Cahn model here considered.  Due
to the noise a solution with center can \emph{nucleate}, i.e., give
rise to new zeroes, during the evolution; in this case it makes no
longer sense to speak of ``the'' center and we simply say \emph{interface
position}.  Such solutions, as we shall see, must be treated carefully
in the statistics.  We note also that an interface may disappear,
either by \emph{exiting} the domain $D$ or by \emph{annihilating}
another interface.

\subsection{A benchmark}
\label{sse:benchmarking}
According to known results \cite[Theorem
  8.1]{funaki:95}\cite{brassesco-demasi-presutti:95} and the
discussion done in \ref{sse:low.noise.resolution}, we expect the
center of $U^n\in\sobh1(D)$, if any, to perform a (discrete) Brownian
Motion, modulo perturbations of order $\Oh(\epsi)$ and the numerical
discretization error.  Therefore, we will declare our numerical scheme
to be acceptable, if we see evidence of this Brownian motion in our
computational results.

Based on this observation, our \emph{benchmarking procedure} consists
in tracking the center $(\Xi^n_\w)_{0\leq n\leq{T/k}}\subset(-1,1)$ of
the computed sample path $(U^n_\w)_{0\leq n\leq{T/k}}$, where $\w$
ranges in a discrete sample space, and we perform statistics on
$(\Xi^n_\w)_{0\leq n\leq{T/k}}$ by averaging over $\w$.  In order to
obtain meaningful statistics, we keep the number of interfaces
constrained to $1$.  Paths that maintain a center up to the final time
are called \emph{admissible sample paths}. If new interfaces are
created, or the center exits from $(-1,1)$, during the computation of
one sample path, it is rejected as not \emph{admissible}, the
computation stopped, and the computation of a new sample path is
started.  The average and variance of the interface position are
computed over the admissible sample paths. The resulting average
interface position, $\EX[\Xi^n]$, and its variance, 
$\Var[\Xi^n]=\EX[(\Xi^n)^2]-(\EX[\Xi^n])^2$, are real valued
functions of (discrete) time $t_n$.

According to \eqref{eqn:funaki.convergence} and the subsequent
remarks, the average position $\EX[\Xi^n]$ must be close to $0$ and its
variance must be close to $c_0\epsi^{1+2\gamma}t_n$, and this is what
we will be after in the next section.

Note that our use of the benchmark is no so strict than the one
usually used in the context numerical schemes for deterministic
nonlinear equations.  In the latter case, an exact solution is usually
readily obtained and the benchmark procedure simply consists in
measuring the error between the two solutions. In the linear
stochastic case, a similar approach can be used with the moments of
the SPDE's solution \cite{allen-novosel-zhang}.  We should stress,
however, that in our case, which is both stochastic and nonlinear, the
analytic knowledge of the exact solution is too scarce, making such
simple benchmarking nearly impossible.

\subsection{Simulations and results}

We run a series of Monte Carlo tests with various combination of
parameters as following:
\begin{algorithmic}
  \STATE $S=0$
  \WHILE{$S\leq S_\mathrm{max}$}
  \FOR{$\w=$clock time} 
    \STATE seed of the random number generator with \w
    \FOR{$\gamma=0.0,\,0.2,\,0.5$}
      \FOR{$\epsilon=0.08,\,0.04,\,0.02,\,0.01$}
        \FOR{$l=6,\,7,\,8,\,9$ \COMMENT{$l$ is the \emph{refinement level}}}
          \STATE let $h=2^{1-l}$ and $k=h^2$
	  \FOR{$n\in\fromto1{T/k}$}
	    \STATE solve \eqref{eqn:newton-fem-sac-matrix} for $U^n_\w$
	    \STATE find $\Xi^n_\w$ such that $U^n_\w(\Xi^n_\w)=0$
	    \IF{$\Xi^n_\w$ exists and is unique}
	      \STATE let $\EX[\Xi^n]_{\mathrm{new}}=(\EX[\Xi^n] S+\Xi^n_\w)/S$
	      \STATE let $\EX[(\Xi^n)^2]_{\mathrm{new}}=(\EX[(\Xi^n)^2] S+\Xi^n_\w)/S$
	    \ELSE
	      \STATE break and skip to the next $l$
	    \ENDIF
	  \ENDFOR
	  \STATE declare sample path successful:
	  \STATE $S=S+1$
	  \STATE $\EX[\Xi^n]=\EX[\Xi^n]_{\mathrm{new}}$
	  \STATE $\EX[(\Xi^n)^2]=\EX[(\Xi^n)^2]_{\mathrm{new}}$
	\ENDFOR
      \ENDFOR
    \ENDFOR
  \ENDFOR
  \ENDWHILE
\end{algorithmic}

The results are reported in Figures \ref{fig:gamma00},
\ref{fig:gamma02} and \ref{fig:gamma05} for the values of
$\gamma=0.0,\,\,0.2,\,0.5$, respectively.  In each figure the
sub-figures (a), (b) and (c), each of which is split into a top and
bottom part, show the graph, as (discrete) functions of $t_n\in[0,T]$,
of the (discrete) average position $\EX[\Xi^n]$ in the top part and
its variance $\var\Xi^n:=\EX[(\Xi^n)^2]-\EX[\Xi^n]^2$ rescaled by
$1/\epsi^{1+2\gamma}T$ in the bottom part for the values of
$\epsi=0.04,\,0.02,\,0.01$, respectively. (For easier visualization we
plot the piecewise linear interpolation of discrete functions.)
Different lines in each of these 3 diagrams correspond to different
values of the mesh refinement level $l$.  The absence of a line means
that the total number of successful sample paths is $0$ and no
statistics are produced.  The diagrams (a), (b), (c) correspond, on each
figure, to the values of $\epsi=0.04,\,0.02,\,0.01$. The plot (d) is
designed so to see how well the scheme captures the
$c_0\epsi^{1+2\gamma}t$ behavior of the variance; each line
corresponds to a chosen but fixed value of $t$, and is a log-log plot
of the variance for the finest refinement level ($l=9$) against
$\epsi$; different lines correspond to different times between $0$ and
$T=20$.  For comparison we plot the line with slope $1+2\gamma$, which
represents the scaling power; lines parallel to this line mean that
the code picks up the right scaling.

We summarize next the observations we have drawn from the computational
results that were described.

\subsubsection{Interface motion}
The motion of the numerical interface, for sufficient transition layer
and noise resolution, has the properties predicted by the analytical
results.  Indeed, starting from a resolved profile centered at zero,
as seen in Figures \ref{fig:gamma00}--\ref{fig:gamma05} the average
position is near zero, whereas the variance, which is expected to be a
linear function of time, behaves in accordance to the expectations, at
least for some initial times; the $\gamma,\epsilon$ dependence of the
diffusion coefficient $c_0\epsi^{1+2\gamma}$ is clearly captured by
the numerics, as seen by the diagram (d) of each of the figures.

\subsubsection{Noise resolution}\label{sse:noise.resolution}
It is well known that for simulating the deterministic Allen--Cahn
equation with any type of mesh/grid-based schemes, the meshsize has to
be smaller than $\epsi$ in order to resolve satisfactorily the
transition layer about the interface.  Intuitively, this is due to the
fact that the transition layer being of width $\Oh(\epsi)$, and the
numerical discretization parameters must be smaller than this width as
to have a proper resolution of this transition layer.

The effect of this is seen in each of Figures \ref{fig:gamma00} to
\ref{fig:gamma05}, as $\epsi$ decreases, the level of refinement has
to be taken bigger and bigger in order to obtain meaningful
calculations.

In the stochastic case, the situation is complicated even more by the
noise.  Indeed, according to Theorem \ref{the:low.noise} and the
discussion in \secref{sse:low.noise.resolution}, the discretization
parameters, in this case $h=\sigma$ and $k=\rho$, must be taken small
enough as to ensure that the noise effects are not lost for small
$\epsi$. Roughly speaking, the discretization parameters must be
small, not only to resolve the interface, but also to resolve the
noise and pick up the diffusion of the Brownian Motion.  This
dependence, which is indicated analytically for $\gamma>3$ by Theorem
\ref{the:low.noise}, is also reflected in our computations for
$\gamma\leq3$.  Indeed, a comparison between Figure \ref{fig:gamma00}
and Figure \ref{fig:gamma02} shows how the level of refinement $l=7$
leads to meaningful results for $\gamma=0.0$ and all values of $\epsi$
(at least for short times), whereas the same refinement level, for the
same values of $\epsi$, but with $\gamma=0.2$ is insufficient.  For
$\gamma=0.5$ in Figure \ref{fig:gamma05} this phenomenon becomes yet
more apparent.

We note that computations with $\gamma>3$ (not shown here) require an
extremely fine mesh, and thus a very small timestep in view of Remark
\ref{obs:convergence-rate}, in order to capture any of the noise
effects.  Otherwise the deterministic solution will be computed.  In
this case, even \emph{choices of $h$ that resolve the interface
satisfactorily are not enough to resolve the noise}.  The lesson we
learn from this, is that the interplay between the noise and the
nonlinearity can be quite delicate, and not obviously predicted from
deterministic considerations, in problems such as the stochastic
Allen--Cahn equation.

\subsubsection{Deterioration of simulations with $t$ big and $\epsi$ small}
\label{sse:deterioration}
As observed in the previous paragraph, the computed variance depends
linearly on time, as expected, but only for some initial time. The
smaller $\epsi$, the shorter this time is.  This is seen in the bottom
part of the sub-figures (a), (b) and (c) by the graph's
earlier or later departure from an initial linear behavior.  The
computations deteriorate faster for small $\epsi$, e.g., $\epsi=0.0$,
than they do for bigger $\epsi$, e.g., $\epsi=0.5$.  The example
with $\epsi=0.2$ shows an intermediate behavior.

\subsection{Interface drift}
To conclude, we add some results of computations, for short times,
with initial value a resolved profile centered away from $0$.  In this
case, the SDE describing the motion of the interface has also a drift
term, which drives the interface towards the closest boundary of the
domain.  This drift is clearly seen for various choices of the
parameters in the top diagrams of Figure \ref{fig:drift}, where we
plot the center's position average against time.  Short times must be
taken, for the statistics to make sense, otherwise solutions with
centers that exit (or nucleate) cease to counterbalance those who stay
in.  Indeed, in the diagram (c), where the noise intensity is quite
strong, after an initial drift towards the boundary, the average
inverts its route and moves away from the boundary.  This is due to
the fact that the statistics become too biased; the SE diagram shows
the samples survival with respect to time (that is, a
discrete-probability space version of the exit-time inverse function).
The samples survival for the top diagrams is 100\% for the statistics
in the top diagrams.  It is worth mentioning that similar
observations, using stochastic ODE's were made by Shardlow
\cite{shardlow:00}.

\newcommand{\XTRIM}{120}
\renewcommand{\picwidth}{.9\textwidth}
\begin{figure}[t]
  \begin{center}
    \subfigure[
      \label{fig:drift.A}
      The average of the interface position (center) for $\gamma=0.0$.
      An upward drift, that is a drift of the interface towards the
      closest boundary point, is visible in this series of
      computations.  The number of the Monte Carlo samples is 1623.
      Only 3 of these paths violate the uniqueness of center condition
      and are therefore excluded from the statistics at the time when
      that happens.
    ]{
    \includegraphics[width=\picwidth, trim = 0 200 125 11, clip]
		    {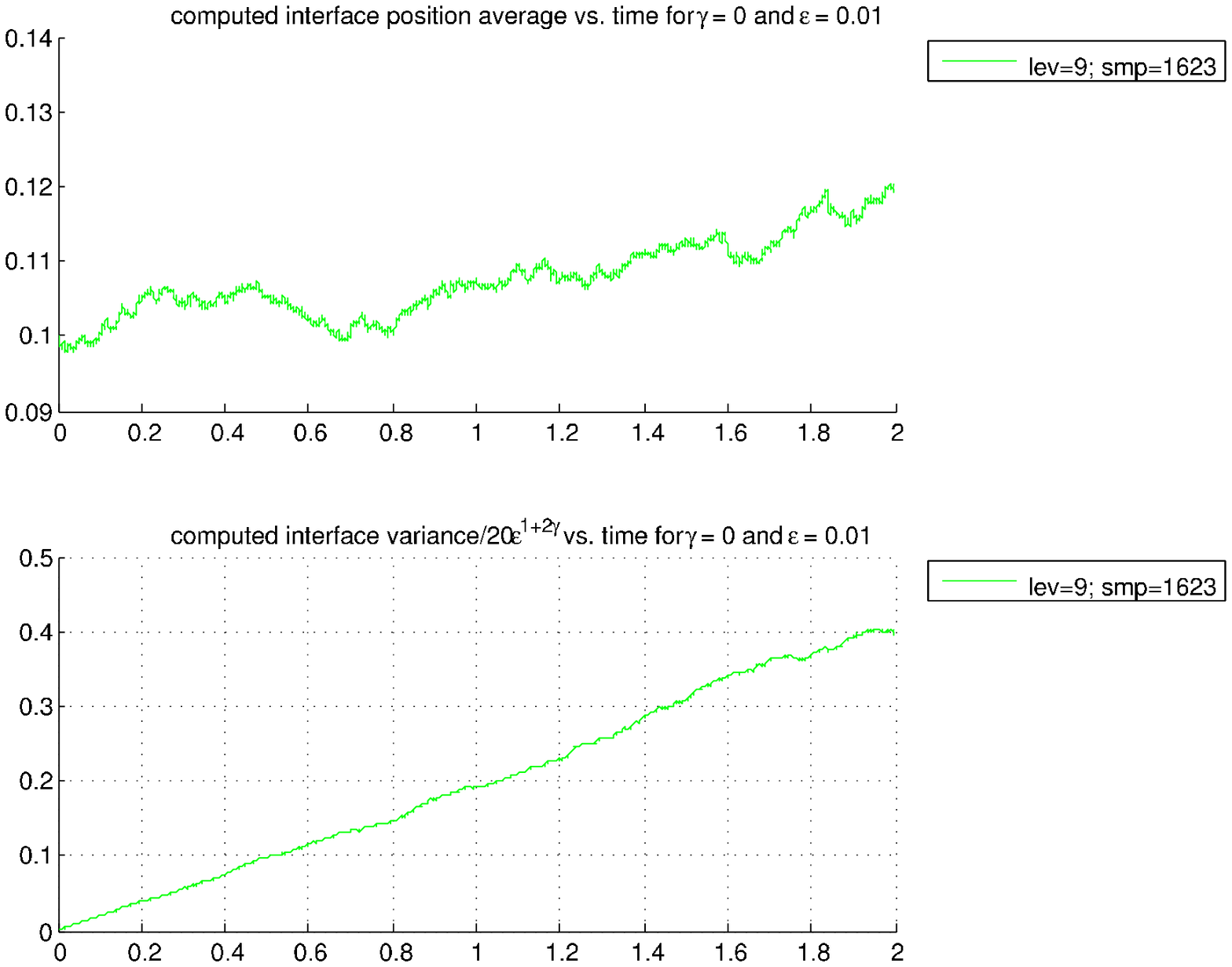}
    }
    \\
    \subfigure[
      The average of the interface position for $\gamma=0.2$.
      Compared with Figure \ref{fig:drift.A}, an outward drift is also
      detectable here but, due to the weaker noise, intensity it is
      slower.  The number of Monte Carlo samples here is 1752, all of
      which make it up to time $T=2$.
    ]{
    \includegraphics[width=\picwidth, trim = 0 200 125 11, clip]
		    {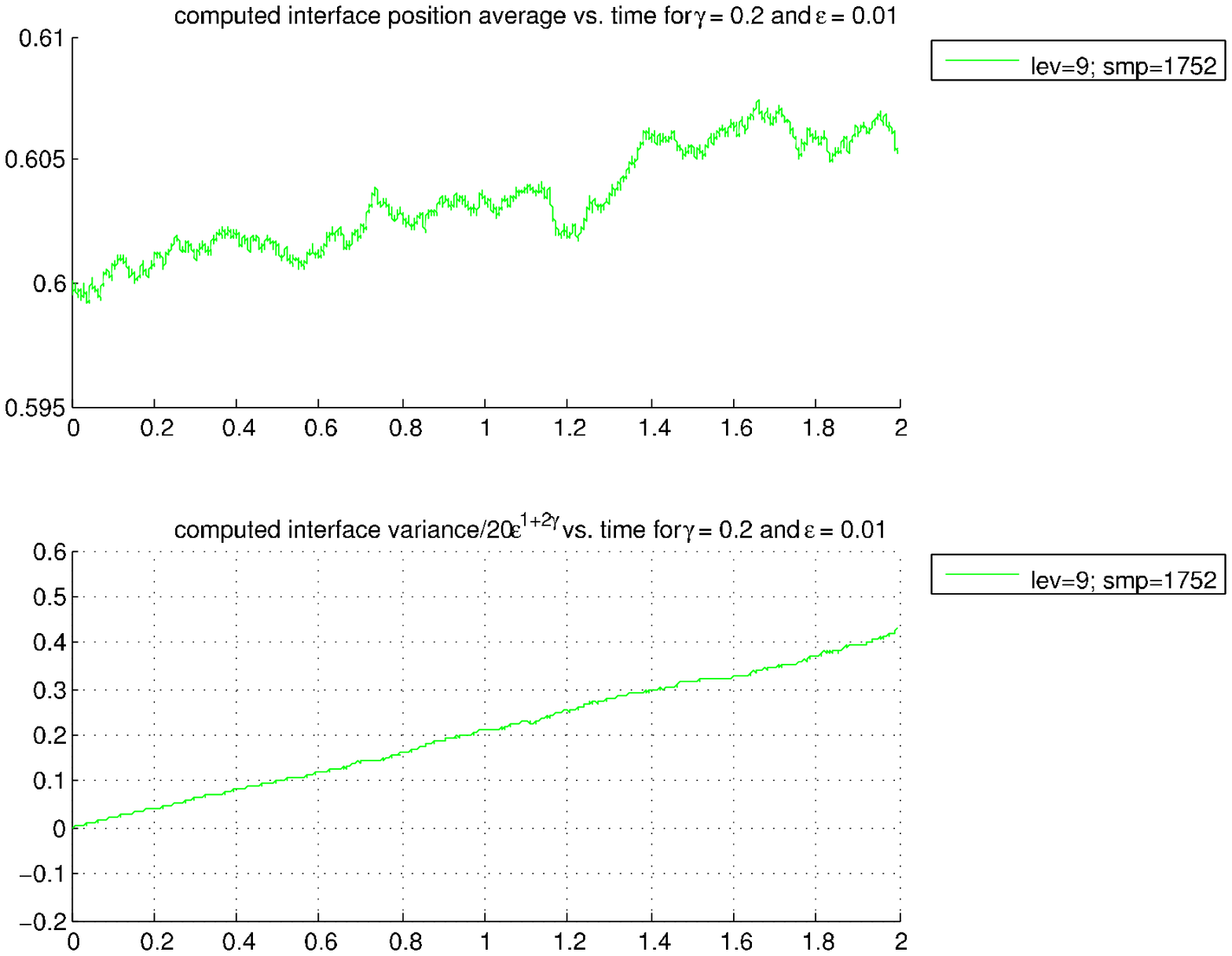}
    }
    \\
    \subfigure[
      On the left, we plot the average of the interface position for
      $\gamma=-0.1$.  In this case, most sample paths do not make it
      with the condition of one interface, or center, up to time
      $T=2$.  This fact is visualized on the right, where the curve
      represents the number of ``admissible'' sample paths (i.e.,
      those that have a center) versus time; as time increases the
      number of the admissible sample paths drops from an initial 2500
      to almost 1500.  This high drop in the number of samples makes
      the statistics unreliable, and explains why the graph on the
      left exhibits no clear outward drift as a result.
    ]{{
      \includegraphics[width=.57\textwidth, trim = 0 200 125 11, clip]
		    {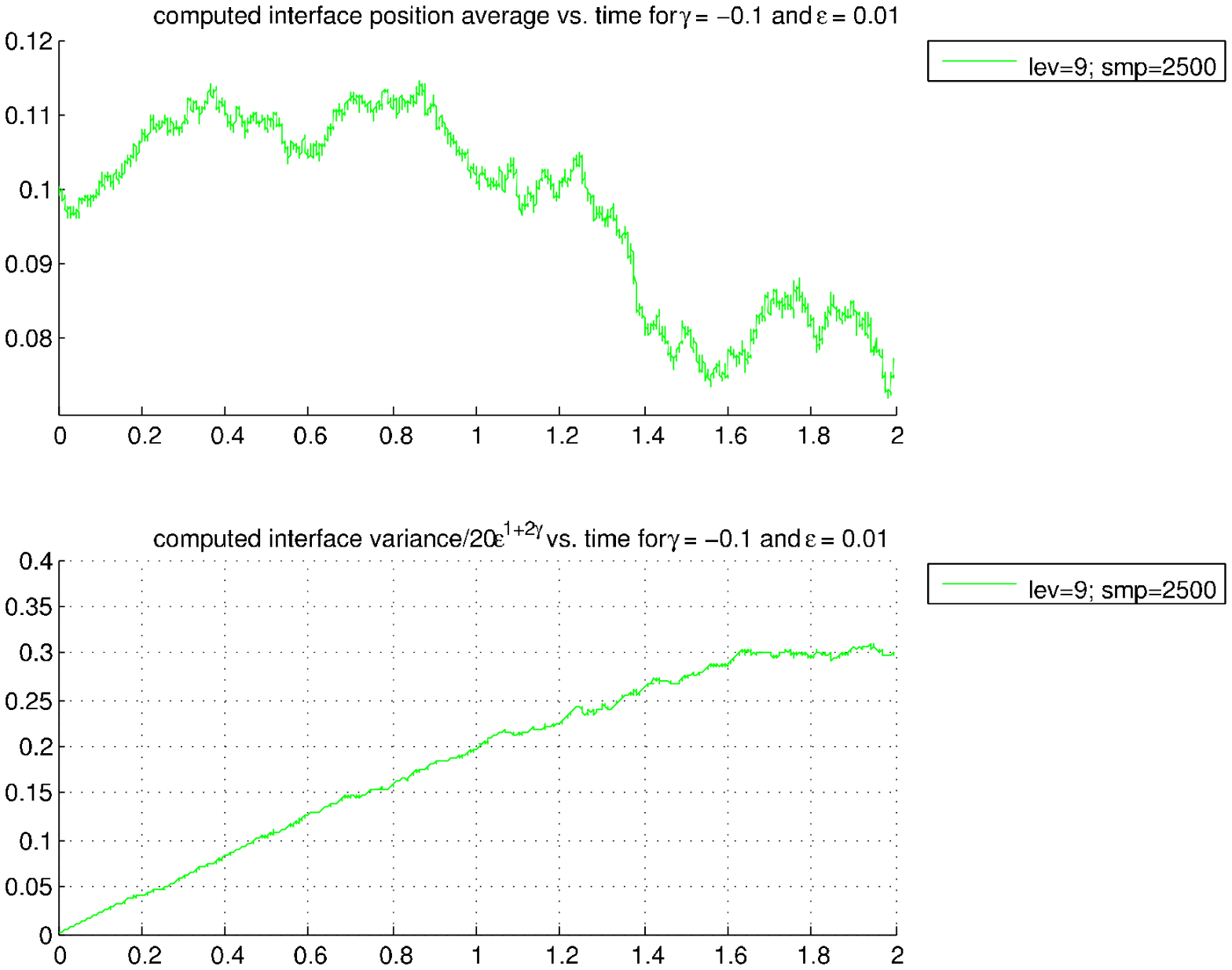}
      \includegraphics[width=.35\textwidth, trim = 0 0 0 0, clip]
		      {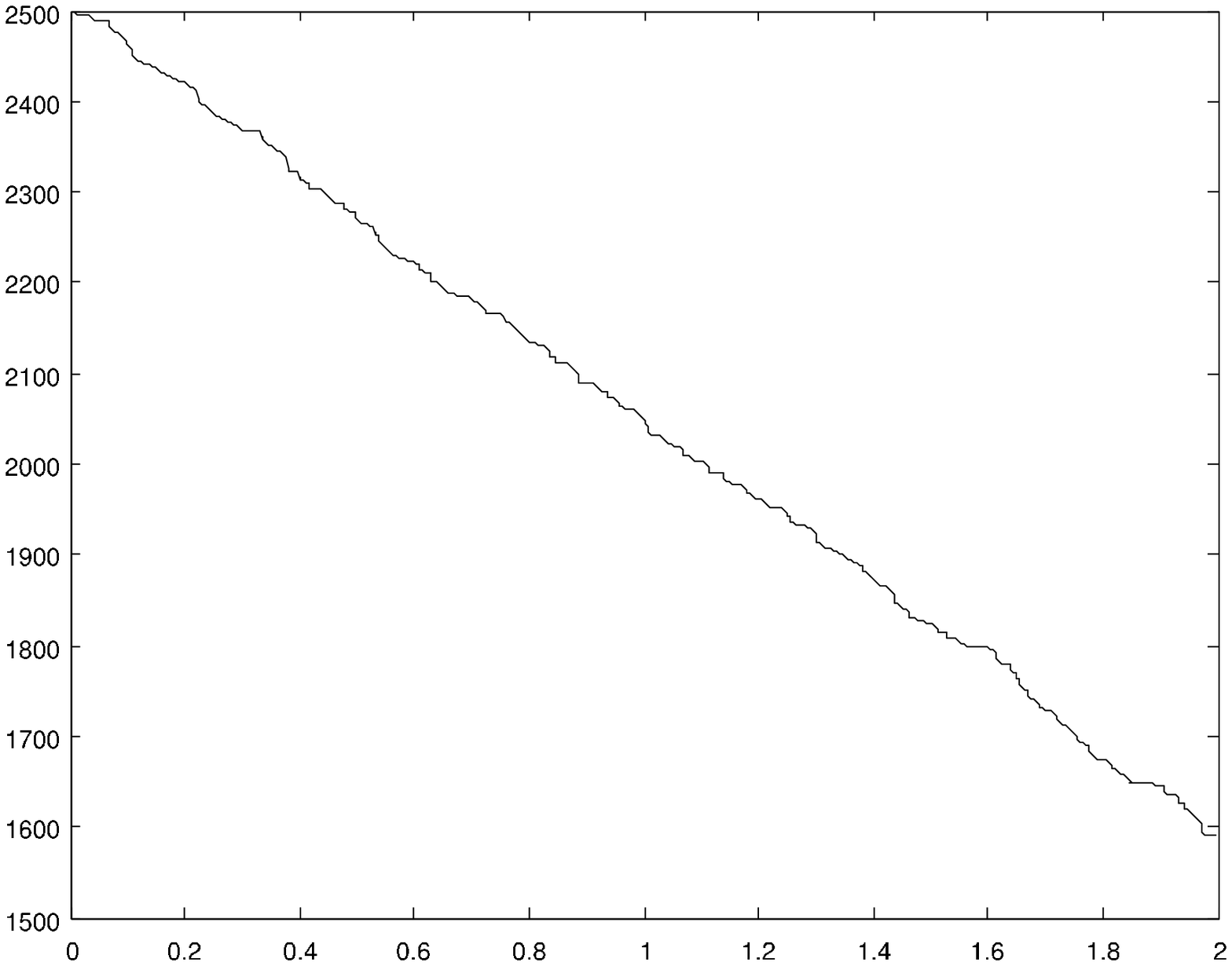}
    }}
  \end{center}
  \caption[Drift computations]{
    \label{fig:drift} Monte Carlo
    computations of the drift effect for the initial interface
    position shifted away from zero.  All these computations are
    performed for $\epsi=0.01$, meshsize $h=2^{-8}$ and timestep
    $k=h^2$.  See the text and Shardlow \cite{shardlow:00}, for
    example, for details about the drift.  
  }
\end{figure}

\subsection*{Acknowledgment}
We would like to thank Georgios Zouraris and Anders Szepessy for
stimulating discussions, and the anonymous referees for their interesting
suggestions that led to the improvement of several aspects of this
paper.

\def\cprime{$'$}
\providecommand{\bysame}{\leavevmode\hbox to3em{\hrulefill}\thinspace}
\providecommand{\MR}{\relax\ifhmode\unskip\space\fi MR }
\providecommand{\MRhref}[2]{%
  \href{http://www.ams.org/mathscinet-getitem?mr=#1}{#2}
}
\providecommand{\href}[2]{#2}

\end{document}